\UseRawInputEncoding

\documentclass[a4paper,hidelinks,11pt]{amsart}
\usepackage{graphicx}
\usepackage{amsmath}
\usepackage{amssymb}

\usepackage{bm} 
\usepackage[T1]{fontenc}
\usepackage{amsfonts}
\usepackage{mathrsfs}
\usepackage{yfonts}
\usepackage{url}

\usepackage[a4paper,margin=3cm]{geometry}

\usepackage{mathtools}
\usepackage{footmisc}

\usepackage{hyperref}
\usepackage{graphicx}
\usepackage{verbatim}
\usepackage{enumerate}	
\usepackage[english]{babel}
\usepackage{ebproof}		
\usepackage{quoting}
\usepackage{array,booktabs}
\usepackage[all]{xy}\SelectTips{cm}{}
\usepackage{pifont}

\usepackage[shortlabels]{enumitem}
\usepackage{units}
\usepackage{xspace}\xspaceaddexceptions{\,}

\usepackage{booktabs,caption}
\usepackage{tikz-cd}
\quotingsetup{font=small}	% Use pdf, png, j

\usepackage{amsthm}

\usepackage{todonotes}

\theoremstyle{definition}   %questo mette le definizione in carattere normale, non corsivetto
\newtheorem{definition}{Definition}[section]

\theoremstyle{plain}    %questo mette i teoremi e i lemmi in corsivetto
\newtheorem{lemma}[definition]{Lemma}
\newtheorem{proposition}[definition]{Proposition}
\newtheorem{theorem}[definition]{Theorem}
\newtheorem{corollary}[definition]{Corollary}

\renewcommand{\phi}{\varphi}

\newcommand\pair[1]{{\langle#1\rangle}}

\newcommand\n{\noindent}

\newcommand\Fm{\mathbf{Fm}}
\newcommand{\logica}{\mathbf}
\newcommand\A{\mathbf{A}}
\newcommand\B{\mathbf{B}}

\newcommand\M{\mathfrak{M}}
\newcommand{\class}{\mathsf}

\newcommand\PL{{\mathcal{P}}_{\textit{\l}}}

\newcommand\WK{{\mathbf{WK}}}
\newcommand\WKt{\WK_{e}}

\makeatletter
\providecommand*{\dashv}{%
	\mathrel{%
		\mathpalette\@dashv\vdash
	}%
}
\newcommand*{\@dashv}[2]{%
	\reflectbox{$\m@th#1#2$}%
}
\makeatother

\makeatletter
\providecommand*{\Dashv}{%
	\mathrel{%
		\mathpalette\@Dashv\vDash
	}%
}
\newcommand*{\@Dashv}[2]{%
	\reflectbox{$\m@th#1#2$}%
}
\makeatother

\newcommand\PWK{\logica{PWK}}
\newcommand\PWKe{\logica{PWK}_{\mathrm{e}}}
\newcommand\MPWK{\ensuremath{\logica{PWK}^{\Box}_{\mathrm{e}}}\xspace}
\newcommand\MPWKglobal{\ensuremath{\logica{PWK}^{g\Box}_{\mathrm{e}}}\xspace}

\newcommand{\IBSL}{\class{IBSL}}
\newcommand{\SIBSL}{\class{SIBSL}}
\newcommand{\BCA}{\class{BCA}}
\newcommand{\MBCABoch}{\class{MBCA_B}}
\newcommand{\MBCAHalld}{\class{MBCA_H}}
\newcommand{\BAO}{\class{BAO}}
\newcommand{\BAOmeet}{\class{BAO_\land}}
\newcommand{\BAOjoin}{\class{BAO_\lor}}

\newcommand{\Jzero}{J_{0}}
\newcommand{\Juno}{J_{1}}
\newcommand{\Jdue}{J_{2}}

\newcommand\ant{\nicefrac12}

\newcommand{\Var}

\renewcommand\B{\logica{B}}
\newcommand\Be{\logica{B}_{\mathrm{e}}}
\newcommand\MBeK{\ensuremath{\logica{B}^{\Box}_{\mathrm{e}}}\xspace}
\newcommand\MBeKglobal{\ensuremath{\logica{B}^{g\Box}_{\mathrm{e}}}\xspace}

\DeclareSymbolFont{symbolsC}{U}{txsyc}{m}{n}
\DeclareMathSymbol{\strictif}{\mathrel}{symbolsC}{74}

\title[]{Modal weak Kleene logics: algebraic semantics and J\'onsson-Tarski representation}
\author[]{S. Bonzio, G. Vergottini and N. Zamperlin}
\date{}

\begin{document}
	
	\begin{abstract}
		In this paper we continue the work started in \cite{bonzio_zamperlin25}	by investigating the algebraic semantics for modal external weak Kleene logics. In order to do so we study the global modal logics associated with the three-valued Kripke models introduced in the previous article, showing how these systems are algebraizable. We proceed by defining their algebraic counterparts, obtaining the classes of modal Bochvar and modal Halld\'{e}n algebras, whose structure theory we investigate via P\l onka sum decomposition. Of particular interested is the interaction of the modal operator with the homorphisms of the P\l onka decomposition. These results allow us to prove representation theorems for these expanded modal algebras in the style of the J\'onsson-Tarski representation for Boolean algebras with operators. As a consequence, we also obtain algebraic semnatics for the local modal logics from which the work originated.
	\end{abstract}
	
\maketitle
\sloppy

Weak Kleene logics are three-valued logics characterized by the presence of an infectious truth value, these include Bochvar's logic ($\mathbf{B}$) and paraconsistent weak Kleene logic ($\mathbf{PWK}$). Despite their logical and semantic interest \cite{CiuniCarrara, Bonzio16SL, beall2016off}, these systems have rather weak algebraic properties. In particular, neither $\mathbf{B}$ nor $\mathbf{PWK}$ is protoalgebraic, and consequently, neither admits an equivalent algebraic semantics in the sense of abstract algebraic logic. This limitation motivates the study of their external versions, obtained by expanding their language with additional unary connectives that establish a stronger connection between logical consequence and algebraic structures.

More precisely, the external versions of these logics are obtained by adding a unary connective $J_2$, which expresses the fact that a formula is classically true. The resulting systems are called Bochvar external logic ($\Be$, introduced by D.A. Bochvar \cite{Bochvar38}) and external paraconsistent weak Kleene logic ($\PWKe$, introduced by S. Halld\'{e}n \cite{Hallden}). Both logics are algebraizable and share the quasivariety $\BCA$ of Bochvar algebras as their equivalent algebraic semantics \cite{SMikBochvar}. The P\l onka-style representation of Bochvar algebras plays an important role in understanding their structure.

The modal extensions of external weak Kleene logics were introduced in \cite{bonzio_zamperlin25}, where the systems $\mathbf{B}_e^\Box$ and $\mathbf{PWK}_e^\Box$ were axiomatized, and it was shown that they are complete with respect to appropriate three-valued Kripke semantics. Both systems are based on the same underlying three-valued operations, but their modal operators are interpreted differently, giving rise respectively to Bochvar-Kripke and Halld\'{e}n-Kripke models. The truth of modal formula $\Box\phi$ has different semantic clauses in the two settings: when $\phi$ is meaningful (i.e. classical) in the current world, the Bochvar modality requires it to be true in every accessible world, whereas the Halld\'{e}n modality requires it simply not to be false in such worlds. In both cases, the infective value of $\Box\phi$ in the current world is only determined by that of $\phi$ in the same worlds, obtaining a form of local infectivity. 

From an algebraic perspective, two questions arise: how does the infective truth value affect the relationship between modal algebras and relational frames, and do the Bochvar and Halld\'{e}n interpretations of modality give rise to different algebraic structures?

One of the main aims of this work is to extend the J\'onsson-Tarski representation theorem for Boolean algebras with operators \cite{jonsson1951boolean, Blue_book} to the context of modal weak Kleene logics. It is the fundamental theorem underlying the algebraization of modal logic and states that every Boolean algebra with operator can be embedded into the complex algebra of a suitable frame built from its ultrafilters. This result establishes a connection between algebraic and relational semantics and is closely related to the construction of the canonical models used in completeness proofs. However, the extension to Bochvar algebras is not immediate, since these algebras admit a P\l onka decomposition into Boolean fibers whose connecting homomorphisms need not preserve the modal operator. Our goal is to obtain analogous representations for the modal algebras associated to $\MBeK$ and $\MPWK$.

We first consider the global consequence relations associated with the two modal logics. Unlike local consequence, global consequence requires that the premises take designated values at every world of a Kripke model and corresponds to a stronger necessitation rule. We show that the resulting global logics, $\MBeKglobal$ and $\MPWKglobal$, are complete with respect to their global Kripke semantics and algebraizable using the same translations as their non-modal counterparts. Their equivalent algebraic semantics are two quasi-varieties of modal expansions of Bochvar algebras, called respectively modal Bochvar algebras ($\MBCABoch$) and modal Halld\'{e}n algebras ($\MBCAHalld$).

A natural question is how the modal operator interacts with the P\l onka decomposition of the underlying Bochvar algebras. We show that $\Box$ preserves each Boolean fiber, although its properties differ in the two classes: in $\MBCABoch$, the operator preserves binary meets within each fiber, whereas every fiber of an algebra in $\MBCAHalld$ is a modal algebra. It is worth to notice that the P\l onka homomorphisms connecting the fibers do not necessarily preserve $\Box$, so these structures cannot, in general, be represented as P\l onka sums in the expanded modal language. Despite this limitation, we show that each fiber, equipped with its modal operator, is isomorphic to a suitable Boolean interval of the bottom fiber, which in turn is equipped with an induced modal operator.

The main representation results of this paper are obtained by introducing two notions of external ultrafilter, starting from the ultrafilters of the bottom Boolean fiber and extending them to the whole Bochvar algebra. These notions reflect the different interpretations of designated values in the Bochvar and Halld\'{e}n semantics. They allow us to associate a relational frame with each algebra and to prove that every algebra in $\MBCABoch$ or $\MBCAHalld$ can be embedded into the corresponding complex algebra of its ultrafilter frame. In this way, we obtain two generalizations of the J\'onsson-Tarski representation theorem. Furthermore, these constructions provide an algebraic characterization of the local modal logics in terms of matrices whose designated sets are the corresponding ultrafilters.

Finally, we analyze the relationship between modal Bochvar and Halld\'{e}n algebras. Although the modal operators have different interpretations, we prove that the quasivarieties $\MBCABoch$ and $\MBCAHalld$ are term-equivalent. More precisely, the modal operator of each class can be defined in terms of the operator of the other, leaving the underlying Bochvar algebra unchanged.

The paper is organized as follows. Section \ref{sct:Preliminaries} presents some concepts related to modal weak Kleene logics, the algebraic properties of Bochvar algebras and P\l onka sums, which are necessary to keep the paper self-contained. In Section \ref{sct:BocHallalgebras}, we introduce modal Bochvar and Halld\'{e}n algebras and analyze their structural properties, with particular attention to the interaction between the modal operator and the P\l onka decomposition. Section \ref{sct:representation} develops the representation theory for these algebras, establishing the corresponding representation theorems via ultrafilter frames. The section concludes by proving the term-equivalence of the two quasivarieties. %Finally, Section \ref{sct:Conclusion} is devoted to concluding remarks.

\section{Preliminaries}\label{sct:Preliminaries}
\subsection{Weak Kleene logics}

Kleene three-valued logics are traditionally divided into two families, depending on the meaning given to the connectives: \emph{strong} Kleene logics -- counting strong Kleene \cite{KleeneBook} and the logic of paradox \cite{Priestfirst} -- and \emph{weak} Kleene logics, namely Bochvar logic and paraconsistent weak Kleene \cite{Bonzio16SL} (sometimes referred to as Halld\'{e}n's logic \cite{Hallden}). All the mentioned four logics are traditionally defined in the language of classical logic. However, the intent of the first developer of the weak Kleene formalism, D.A. Bochvar, was to work within an enriched language allowing to express all classical ``two-valued'' formulas -- which he referred to as \emph{external formulas} -- beside the genuinely ``three-valued'' ones. 

The result of this choice is the language $\langle \neg,\lor,\Jdue,0,1\rangle$ of type $\pair{1,2,1,0,0}$, which is obtained by enriching the classical language by an additional unary connective $\Jdue$ (and the constants $0,1$), where the formula $\Jdue\varphi$ is to be read as ``$\varphi$ is true''. This language will be referred to as \emph{external language}, while its $\Jdue$-reduct will be called \emph{internal}. 

The algebraic interpretation of the language is given via the three-element algebra $\WKt=\pair{\{0,1,\ant\}, \neg,\lor,\Jdue,0,1}$ displayed in Figure \ref{fig:WKe}.
\begin{figure}[h]
	\begin{center}\renewcommand{\arraystretch}{1.25}
		\begin{tabular}{>{$}c<{$}|>{$}c<{$}}
			& \lnot  \\[.2ex]
			\hline
			1 & 0 \\
			\ant & \ant \\
			0 & 1 \\
		\end{tabular}
		\qquad
		\begin{tabular}{>{$}c<{$}|>{$}c<{$}>{$}c<{$}>{$}c<{$}}
			\lor & 0 & \ant & 1 \\[.2ex]
			\hline
			0 & 0 & \ant & 1 \\
			\ant & \ant & \ant & \ant \\          
			1 & 1 & \ant & 1
		\end{tabular}
		\qquad
		\begin{tabular}{>{$}c<{$}|>{$}c<{$}}
		& J_{_2} \\[.2ex]
			\hline
			1 & 1 \\
			\ant & 0 \\
			0 & 0 \\
		\end{tabular}
	\end{center}
	\caption{The algebra $\WKt$.}\label{fig:WKe}
\end{figure}

\n The third value $\ant$ is traditionally read as ``meaningless'' (see e.g. \cite{ferguson2017meaning} and \cite{SmuczMeaningless}) due to its infectious behavior. The connectives of conjunction and material implication are defined in the standard way and have the expected tables:

\begin{figure}[h]
	\begin{center}\renewcommand{\arraystretch}{1.25}
		\begin{tabular}{>{$}c<{$}|>{$}c<{$}>{$}c<{$}>{$}c<{$}}
			\land & 0 & \ant & 1 \\[.2ex]
			\hline
			0 & 0 & \ant & 0 \\
			\ant & \ant & \ant & \ant \\          
			1 & 0 & \ant & 1
		\end{tabular}
		\qquad
		\begin{tabular}{>{$}c<{$}|>{$}c<{$}>{$}c<{$}>{$}c<{$}}
			\to & 0 & \ant & 1 \\[.2ex]
			\hline
			0 & 1 & \ant & 1 \\
			\ant & \ant & \ant & \ant \\          
			1 & 0 & \ant & 1
		\end{tabular}
	\end{center}
\end{figure} 

Via $\Jdue$ we can define other external connectives: $\Jzero\varphi\coloneqq \Jdue\neg\varphi$ expressing the falsity of a formula and $\Juno\varphi\coloneqq \neg(\Jdue\varphi\lor\Jdue\neg\varphi)$ expressing its non-classicality. Moreover, $\varphi\equiv\psi\coloneqq\displaystyle(\Jdue\phi\leftrightarrow\Jdue\psi)\land(\Jzero\phi\leftrightarrow\Jzero\psi)$ is the proper logical equivalence for the external language. Their interpretation in $\WKt$ is displayed in the following tables:
\vspace{5pt}
\begin{center}
	\begin{tabular}{>{$}c<{$}|>{$}c<{$}}
		& \Jzero\\[.2ex]
		\hline
		1 & 0 \\
		\ant & 0 \\
		0 & 1 \\
	\end{tabular}
	\qquad
	\begin{tabular}{>{$}c<{$}|>{$}c<{$}}
		& \Juno\\[.2ex]
		\hline
		1 & 0 \\
		\ant & 1 \\
		0 & 0 \\
	\end{tabular}
	\qquad
	\begin{tabular}{>{$}c<{$}|>{$}c<{$}>{$}c<{$}>{$}c<{$}}
		\equiv & 0 & \ant & 1 \\[.2ex]
		\hline
		0 & 1 & 0 & 0 \\
		\ant & 0 & 1 & 0 \\          
		1 & 0 & 0 & 1
	\end{tabular}
\end{center}
\vspace{5pt}

According to the above interpretation some formulas are evaluated \emph{exclusively} into $\{0,1\}$  (which is the universe of a Boolean subalgebra of $\WKt$), by any homomorphism $h\colon\Fm\to\WKt$: these are called \emph{external} formulas.  

We recall that a logic $\mathsf{L}$ is induced by a matrix $\pair{\mathbf{M},F}$ (of the same language) when
$$\Gamma\vdash_\mathsf{L}\varphi \text{ iff for all homomorphisms } h\colon\Fm\to\mathbf{M},$$ $$h[\Gamma]\subseteq F \text{ implies } h(\varphi)\in F.$$

\begin{definition}\label{def.: logiche esterne}
	Bochvar \emph{external} logic $\Be$ is the logic induced by the matrix $\pair{\WKt,\{1\}}$. PWK \emph{external} logic $\PWKe$ is the logic induced by the matrix $\pair{\WKt,\{1,\ant\}}$.
\end{definition}
%\medskip
%\noindent

Hilbert-style axiomatizations of the two logics are given by Finn-Grigolia \cite{FinnGrigolia} and Segerberg \cite{Segerberg65}, respectively. An important feature of these axiomatizations that will play a crucial role in the following is that a precise fragment of external weak Kleene logics behaves exactly as classical logic. In order to express that we need the notion of external formula.

\begin{definition}
	An occurrence of a variable $x$ in a formula $\phi$ is \emph{covered} if it falls under the scope of any $\Jdue$ operator. A formula $\phi$ is called \emph{external} if all its variables are covered.
\end{definition}

\n The next proposition is immediate by observing that the axiomatizations of external weak Kleene logics contain all the axioms of classical logic when all formulae are taken external and that \textit{modus pones} works for such formulae.

\begin{proposition}\label{prop: classical_fragment_eWKL}
	For a formula $\phi$ define $c(\phi)$ as the formula obtained by removing every occurrence of $\Jdue$. Let $\Gamma,\phi$ be external formulae, then $\Gamma\vdash_{\logica{L}}\phi\Leftrightarrow c[\Gamma]\vdash_{\logica{CL}}c(\phi)$, for $\logica{L}\in\{\Be,\PWKe\}$.
\end{proposition}

\subsection{Modal weak Kleene logics}

In \cite{bonzio_zamperlin25}, possible world semantics are introduced for modal weak Kleene logics. The starting point is a three-valued Kripke model, which is based on a standard Kripke frame, now equipped with a valuation function that assigns to formulae elements of $\WKt$. All the unproved results stated in this section can be found in the mentioned paper. 

In the following, if not otherwise stated, we will be working with the language $\pair{\lor,\neg,\Jdue,\Box,0,1}$ of type $\langle 2,1,1,1,0,0\rangle$. The absolutely free algebra in this language (over denumerable generators) will be denoted by $\Fm$.

\begin{definition}\label{def: Kripke model}
	A \emph{weak three-valued Kripke model} $\M$ is a structure $\pair{W,R, v}$ such that $W$ is a non-empty set (of possible worlds); $R$ is a binary relation over $W$; $v\colon W\times \Fm\to\WKt$ is a homomorphism for the $\Box$-free fragment of the language. 
\end{definition} 

\n Remember the given a class of Kripke models two different notions of logical consequence are commonly definable:
\begin{definition}\label{def: conseguenza logica Bochvar}
	For the class of Bochvar-Kripke models, the \textit{local} consequence relation $\Gamma \vDash^{l}_{\B^{\square}_\mathrm{e}} \varphi$ holds iff for all models $\pair{W,R, v}$ and all $w\in W$, if $v(w,\gamma)=1, \forall\gamma\in\Gamma$, then $v(w,\varphi)=1$. The \textit{global} consequence relation $\Gamma \vDash^{g}_{\B^{\square}_\mathrm{e}} \varphi$ holds iff for all models $\pair{W,R, v}$, if $v(w,\gamma)=1, \forall\gamma\in\Gamma$ and $\forall w\in W$, then $v(w,\varphi)=1,\forall w\in W$.
\end{definition}

In order to adapt the interpretation of the modality $\Box$ to the semantics of Bochvar logic the following class of models has been selected:

\begin{definition}\label{def: semantica Box}
	A \emph{Bochvar-Kripke model} is a weak three-valued Kripke model $\pair{W,R, v}$ such that $v$ evaluates formulae of the form $\Box\varphi$ according to:
	\begin{itemize}
		\item[(1)] $v(w,\Box\varphi)= 1$ iff $v(w,\varphi)\neq\ant$ and $v(s,\varphi)= 1$ for every $s\in W$ such that $wRs$. 
		\item[(0)] $v(w,\Box\varphi)= 0$ iff $v(w,\varphi)\neq\ant$ and there exists $s\in W$ such that $wRs$ and $v(s,\varphi)\neq 1$. 
		\item[($\ant$)] $v(w,\Box\varphi)= \ant$ iff $v(w,\varphi)=\ant$. 
	\end{itemize}
\end{definition}

The local modal logic of Bochvar-Kripke models has been axiomatized by the Hilbert-style system $\B^\square_\mathrm{e}$, which has expected contains a rule of local necessitation:
\begin{itemize}
	\item[] \textbf{Axioms}
	\item the axioms of $\B_\mathrm{e}$ (see \cite{FinnGrigolia}); 
	\vspace{3pt}
	\item[(B1)]  $\square(J_2\varphi\to J_2\psi)\to (\square J_2\varphi\to\square J_2\psi)$;
	\vspace{3pt}
	\item[(B2)]  $\neg\Juno\varphi\leftrightarrow \neg\Juno\square \varphi$; 
	\vspace{3pt}
	\item[(B3)] $\Jdue\square\varphi\to\square\Jdue \varphi$;
	\vspace{3pt}
	\item[(B4)] $\Jzero\square\varphi\to\neg\square\Jzero\neg \varphi$.
	\item[] \textbf{Rules}
	%\item (MP): $\varphi,\varphi\to\psi\vdash\psi$;
	\item the rules of $\B_\mathrm{e}$;
	\vspace{3pt}
	\item[(N)] if $\vdash \varphi $ then $\vdash \square\varphi$.
\end{itemize}

\n It is interesting to notice that this logic coincides with the minimal normal modal logic \textbf{K} when restricted to external formulae, by the presence of the rule of necessitation together with the following lemma.
\begin{lemma}\label{lemma: K per Bochvar su esterne}
	For $\alpha,\beta$ external formulae, the following is a theorem of $\B^{\square}_\mathrm{e}$:
	\begin{itemize}
		\item[\emph{(BK)}] $\square(\alpha\to\beta)\to (\square\alpha\to\square\beta)$
	\end{itemize}    
\end{lemma}

\n We obtain the adaptation of Proposition \ref{prop: classical_fragment_eWKL} for the modal setting.
\begin{lemma}\label{lem: modal_fragment_eWKL}
	For $\Gamma,\phi$ external formulae, $\Gamma\vdash_{\MBeK}\phi\Leftrightarrow c[\Gamma]\vdash_{\logica{K}}c(\phi)$, using the definition of $c(\phi)$ from Proposition \ref{prop: classical_fragment_eWKL}.
\end{lemma}

The completeness proof makes use of a canonical model construction based on the following steps.

\begin{definition}\label{def: modello canonico}
	The \emph{canonical model} for $\B^\square_\mathrm{e}$ is a model $\M=\pair{\mathcal{W},\mathcal{R}, v}$ where $\mathcal{W}$ is the set of all maximally consistent sets of formulae, the canonical relation $\mathcal{R}$ is defined, for every $w,s\in \mathcal{W}$ as
	$$w\mathcal{R} s \text{ iff }\forall\varphi\in Fm \text{ s.t. } \Box\varphi\in w \text{ then } \varphi\in s. $$
	The valuation $v$ is defined as follows, for every $w\in\mathcal{W}$ and variable $x$: 
	\begin{itemize}
		\item $v(w,x)=1$ if and only if $x\in w$;
		\vspace{2pt}
		\item $v(w,x)=0$ if and only if $\neg x\in w$;
		\vspace{2pt}
		\item $v(w,x)=\ant$ if and only if $\Juno x\in w$.
	\end{itemize}
\end{definition}

\begin{lemma}[Lindenbaum's Lemma]\label{lemma: Lindenbaum}
	Let $\Gamma$ be a consistent set of formulae such that $\Gamma\nvdash_{\B^\square_\mathrm{e}}\varphi$, for some $\varphi\in Fm$, then there exists a maximally consistent set of formulae $w$ such that $\Gamma\subseteq w$ and such that $\varphi\not\in w$.
\end{lemma}

\begin{lemma}[Existence Lemma]\label{lemma: di esistenza}
	For every maximally consistent set of formulae $w\in\mathcal{W}$ such that $\diamond\varphi\in w$ then $\varphi$ is meaningful in $w$ (that is $x\in w$ or $\neg x\in w$, for every variable $x\in \varphi$ which does not fall under the scope of $\Jdue$) and there exists a maximally consistent set $s\in\mathcal{W}$ such that $w\mathcal{R}s$ and either $\varphi\in s$ or $\Juno\varphi\in s$. 
\end{lemma}

\begin{lemma}[Truth Lemma]\label{lemma: truth lemma}
	Let $\M=\pair{\mathcal{W},\mathcal{R}, v}$ be the canonical model for $\B^\square_\mathrm{e}$. Then, for every formula $\varphi\in Fm$ and every $w\in\mathcal{W}$, the following hold: 
	
	\begin{enumerate}
		\item $v(w,\varphi)=1$ if and only if $\varphi\in w$;
		\vspace{2pt}
		\item $v(w,\varphi)=0$ if and only if $\neg \varphi\in w$;
		\vspace{2pt}
		\item $v(w,\varphi)=\ant$ if and only if $\Juno\varphi\in w$.
	\end{enumerate}
\end{lemma}

\begin{theorem}[Completeness]\label{th: completezza Bochvar modale}
	$\Gamma\vdash_{\B^\square_\mathrm{e}}\varphi$ if and only if $ \Gamma\vDash^l_{\B^\square_\mathrm{e}}\varphi$.
\end{theorem}

Among the weak three-valued Kripke models of Definition \ref{def: Kripke model} we can choose a different class which better seems to return the correct reading of modality in paraconsistent weak Kleene logic.

\begin{definition}\label{def: semantica Box per PWK}
	A \emph{Halld\'{e}n-Kripke} model is a weak three-valued Kripke model $\pair{W,R, v}$ such that $v$ evaluates formulae of the form $\Box\varphi$ according to:
	\begin{itemize}
		\item[(1)] $v(w,\Box\varphi)= 1$ iff $v(w,\varphi)\neq\ant$ and $v(s,\varphi)\neq0$ for every $s\in W$ such that $wRs$. 
		\item[(0)] $v(w,\Box\varphi)= 0$ iff $v(w,\varphi)\neq\ant$ and there exists $s\in W$ such that $wRs$ and $v(s,\varphi)=0$. 
		\item[($\ant$)] $v(w,\Box\varphi)= \ant$ iff $v(w,\varphi)=\ant$. 
	\end{itemize}
\end{definition}

\n As before, we consider two relations of logical consequence.

\begin{definition}\label{def: conseguenza logica pwk}
	For the class of Halld\'{e}n-Kripke models, the \textit{local} consequence relation $\Gamma \vDash^{l}_{\MPWK} \varphi$ holds iff for all models $\pair{W,R, v}$ and all $w\in W$, if $v(w,\gamma)\neq 0, \forall\gamma\in\Gamma$, then $v(w,\varphi)\neq 0$. The \textit{global} consequence relation $\Gamma \vDash^{g}_{\MPWK} \varphi$ holds iff for all models $\pair{W,R, v}$, if $v(w,\gamma)\neq 0, \forall\gamma\in\Gamma$ and $\forall w\in W$, then $v(w,\varphi)\neq 0,\forall w\in W$.
\end{definition}

The following axiomatic system is the local modal logic of the class of Halld\'{e}n-Kripke models:

\begin{itemize}
	\item[] \textbf{Axioms} 
	\item the axioms for $\PWKe$ (see \cite{Segerberg65}); 
	\vspace{3pt}
	\item[(P1)]  $\square(J_2\varphi\to J_2\psi)\to (\square J_2\varphi\to\square J_2\psi)$;
	\vspace{3pt}
	\item[(P2)]  $\neg\Juno\varphi\leftrightarrow \neg\Juno\square \varphi$.
	\vspace{3pt}
	\item[(P3)]  $\Box\varphi\leftrightarrow \Box\neg\Jzero \varphi$;
	\item[] \textbf{Rules} 
	\item the rules for $\PWKe$; 
	\vspace{3pt}
	\item[(N)] if $\vdash \varphi $ then $\vdash \square\varphi$.
\end{itemize}

\n This logic enjoys the same property of $\MBeK$ described in Lemma \ref{lemma: K per Bochvar su esterne} when it comes to external formulae, therefore we expand Proposition \ref{prop: classical_fragment_eWKL} to $\MPWK$ as well.

\begin{lemma}
	For $\Gamma,\phi$ external formulae, $\Gamma\vdash_{\MPWK}\phi\Leftrightarrow c[\Gamma]\vdash_{\logica{K}}c(\phi)$, using the definition of $c(\phi)$ from Proposition \ref{prop: classical_fragment_eWKL}.
\end{lemma}

In this paraconsistent setting, the completeness proof substitutes the notion of consistency with that of non-triviality ($\Gamma\subset Fm$ is non-trivial when there is no formula $\varphi$ for which $\Gamma\vdash\Jdue\varphi$ and $\Gamma\vdash\neg\Jdue\varphi$)\footnote{Some definitions and lemmas have been corrected since \cite{bonzio_zamperlin25}. The respective proofs still follow the same strategies employed in the original paper.}.

\begin{definition}\label{def: modello canonico MPWK}
The \emph{canonical model} for \MPWK is a weak three-valued Kripke model $\M=\pair{\mathcal{W},\mathcal{R}, v}$ where $\mathcal{W}$ is the set of all maximally non-trivial sets of formulae, $\mathcal{R}$ is the canonical relation for \MPWK defined as, for every $w,s\in \mathcal{W}$ as
$$w\mathcal{R} s \text{ iff }\forall\varphi\in Fm \text{ s.t. } \Box\varphi\in w \text{ and }\Juno\varphi\notin w, \text{ then } \varphi\in s,$$	
while $v$ is defined as follows: 
\begin{itemize}
	\item $v(w,x)=1$ if and only if $\Jdue x\in w$;
	\vspace{2pt}
	\item $v(w,x)=0$ if and only if $\Jzero x\in w$;
	\vspace{2pt}
	\item $v(w,x)=\ant$ if and only if $\Juno x\in w$.
\end{itemize}
\end{definition}

\begin{lemma}[Lindenbaum's lemma]\label{lemma: Lindenbaum per PWK}
	Let $\Gamma$ be a non-trivial set of formulae such that $\Gamma\nvdash_{\MPWK}\varphi$, for some $\varphi\in Fm$, then there exists a complete and non-trivial set of formulae $w$ such that $\Gamma\subseteq w$ and $\neg\varphi\in w$.
\end{lemma}

\begin{lemma}[Existence lemma]\label{lemma: di esistenza per MPWK}
	For every maximally non-trivial set of formulae $w$ such that $\diamond\varphi\in w$ and $\varphi$ is meaningful in $w$, there exists a maximally non-trivial set $s\in\mathcal{W}$ such that $w\mathcal{R}s$ and $\Jdue\varphi\in s$. 
\end{lemma}

\begin{theorem}[Completeness]\label{th: completezza MPWK}
	$\Gamma\vdash_{\MPWK}\varphi$ if and only if $ \Gamma\vDash^{l}_{\MPWK}\varphi$.
\end{theorem}

\subsection{Bochvar algebras and P\l onka sums}

Our discussion over Bochvar algebras will be mostly performed using their representation as P\l onka sums, therefore we start by introducing this technique. A P\l onka sum (\cite{Plo67}, see also \cite{Plo68}, \cite{Plonka69} and the survey in \cite{RomanowskaSurvey92}) is a pasting of similar algebras arranged in a special enriched semilattice.

\begin{definition}
	A \textit{semilattice direct system} of algebras of type $\tau$ is a triple $\mathbb{A} = \pair{\{\A_i\}_{i\in I}, \pair{I,\leq}, \{p_{ij}|i,j\in I, i\leq j\}}$ consisting of:
	\begin{itemize}
		\item a semilattice $\pair{I,\leq}$ with join $\lor$;
		
		\item a family $\{\A_i\}_{i\in I}$ of $\tau$-algebras with pairwise disjoint universes;
		
		\item a homomorphism $p_{ij}:\A_i\to\A_j$, for every $i,j\in I$ such that $i \leq j$. Moreover, if $i\leq j,j\leq k, p_{ik}=p_{jk}\circ p_{ij}$, and $p_{ii}=id_{A_i}$.
	\end{itemize}
	
	\n When $\tau$ contains constants, $\pair{I,\leq}$ is a semilattice with bottom element $0$.
\end{definition}

\begin{definition}\label{def.: plonka_sum}
	For a semilattice direct system $\mathbb{A} = \pair{\{\A_i\}_{i\in I}, \pair{I,\leq}, \{p_{ij}|i,j\in I, i\leq j\}}$, the \textit{P\l onka sum} $\PL(\mathbb{A})$ over $\mathbb{A}$ is the algebra, also of type $\tau$, such that:
	\begin{itemize}
		\item its universe is $\bigcup_{i\in I} A_i$ (each $\A_i$ is called a \textit{fiber} of the sum);
		
		\item for every $n$-ary operation $g\in\tau$, and $a_1\in A_{i_1}, \dots, a_n\in A_{i_n}$;
		
		$$g^{\PL(\mathbb{A})}(a_1,\dots, a_n):=g^{A_j}(p_{i_1j}(a_1),\dots, p_{i_nj}(a_n))$$
		
		where $j=i_1 \lor \dots$ $\lor i_n$;
		
		\item with constants in the type: $c^{\PL(\mathbb{A})}=c^{A_0}$.
	\end{itemize}
\end{definition}

\n The way in which operations are computed once we paste a semilattice direct system into a P\l onka sum is such that if all the arguments of the operation belong to the same fiber, the output is the same as if it would have been computed in the original algebra that constitutes the fiber; if the elements belong to different fibers, they are mapped by the homomorphisms $p_{ij}$ to the fiber which is their least upper bound and the operation is computed there on the respective images.

P\l onka sums resulted greatly helpful in the study of \textit{regular varieties}, those varieties whose valid equations are only of the form $\alpha\approx\beta$ with $Var(\alpha)=Var(\beta)$ (called \textit{regular equations}).

\begin{theorem}[\cite{Plo67}, Theorem 1]\label{th.: plonka_valid_equations}
	If $\mathbb{A}$ is a semilattice direct system containing at least two fibers, the equations valid in $\PL(\mathbb{A})$ are precisely the regular equations that hold in all algebras of $\mathbb{A}$.
\end{theorem}

\n The \textit{regularization} $\mathcal{R}(\class{V})$ of a variety $\class{V}$ is the class of algebras whose equational theory is composed by all and only the regular equations that are valid in $\class{V}$. A variety $\class{V}$ is called \textit{strongly irregular} if there is a term $\phi(x,y)$ in exactly two variables s.t. $\class{V}\vDash\varphi(x,y)\approx x$; e.g., the variety of lattices satisfies the equation $x\land(x\lor y)\approx x$ and as such is strongly irregular, as are all the varieties of algebras with a lattice reduct, e.g., the variety $\class{BA}$ of Boolean algebras.

\begin{theorem}[\cite{RomanowskaSurvey92}, Theorem 7.1]\label{th.: plonka_regularization}
	Let $\class{V}$ be a strongly irregular variety. Then $\mathbf{A}\in\mathcal{R}(\class{V})$ iff $\mathbf{A}$ is decomposable as a P\l onka sum over a semilattice direct system of algebras in $\class{V}$.
\end{theorem}

The regularization of Boolean algebras $\class{BA}$, formulated in the language $\pair{\lor,\land,\neg,0,1}$ of type $\langle 2,2,1,0,0\rangle$, is the class of involutive bisemilattices $\IBSL$ (see \cite[def. 2.4.1]{Bonziobook}), axiomatised as follows:
\begin{enumerate}
	\item $x\vee x\approx x$; 
	\item$x\lor y \approx y \lor x$;
	\item $(x\lor y)\lor z\approx  x\lor(y\lor z)$;
	\item $\neg\neg x\approx x$;
	\item $x\land y\approx \neg(\neg x\vee\neg y)$;
	\item $x\land(\neg x\vee y)\approx x\land y$;
	\item $0\vee x\approx x$;
	\item $1\approx \neg 0$;
\end{enumerate}

\noindent Using the above characterization of $\IBSL=\mathcal{R}(\class{BA})$, it is easy to visualize an involutive bisemilattice, since it is just a semilattice direct system of Boolean algebras. In particular the algebra $\WK$ described in the introduction is an involutive bisemilattice, moreover it is the generator of the entire variety: $\IBSL=\mathbb{V}(\WK)$ (\cite{Bonzio16SL}).

In the case of a P\l onka sum whose fibers are only Boolean algebras, we are going to employ the following abbreviations for any index $i\in I$:

\begin{itemize}
	\item $1_i:=x\lor\neg x$ and $0_i:=x\land\neg x$ for arbitrary $x\in A_i$;
	\item for $x\in A_i,1_x:=1_i$ and $0_x:=0_i$;
	\item if $a\in A_i,A_a:=A_i$.
\end{itemize}

\n $1_i$ represents the top element of the Boolean fiber $\A_i$, in this sense we will also call it the local top, to distinguish it from the algebraic constant $1$ (which is the local top of the bottom fiber $\A_0$). By the same reasoning, $0_i$ represents the bottom element of $\A_i$ and we will call it the local bottom, to distinguish it from $0$. That these writings are not ambiguous relies on the definition of operations in a P\l onka sum and the fact that all fibers are assumed to be Boolean algebras.

Returning to the main topic of weak Kleene logics, the reason behind this algebraic detour is that the algebraic counterpart of external weak Kleene logic, Bochvar algebras, are a particular class of expanded involutive bisemilattices. Bochvar algebras were originally axiomatized by \cite{FinnGrigolia}. As before, we employ the abbreviations $\Jzero\varphi\coloneqq \Jdue\neg\varphi$ and $\Juno\varphi\coloneqq \neg(\Jdue\varphi\lor\Jdue\neg\varphi)$.

\begin{definition}\label{def. bochvar_algebra}
	A \textit{Bochvar algebra} is an structure $\A = \langle A,\neg,\lor,\land,\Jdue,0,1\rangle$ of type $\langle 1,2,2,1,0,0\rangle$ satisfying the following equations and quasiequations:
	\begin{enumerate}
		\item $x\vee x\thickapprox x$; \label{BCA:1}
		\item$x\lor y \thickapprox y \lor x$; \label{BCA:2}
		\item $(x\lor y)\lor z\thickapprox  x\lor(y\lor z)$; \label{BCA:3}
		\item $x\land(y\lor z)\thickapprox(x\land y)\lor(x\land z)$; \label{BCA:4}
		\item $\neg\neg x\thickapprox x$; \label{BCA:5}
		\item $\neg 1\thickapprox 0$; \label{BCA:6}
		\item $\neg( x\lor y)\thickapprox \neg x\land\neg y$; \label{BCA:7}
		\item $0\vee x\thickapprox x$; \label{BCA:8}
		\item $\Jzero\Jdue x\thickapprox \neg\Jdue x$; \label{BCA:10}
		\item $\Jdue x \thickapprox\neg(\Jzero x\vee \Juno x)$; \label{BCA:13}
		\item $\Jdue x\vee\neg \Jdue x \thickapprox 1$; \label{BCA:14}
		\item $J_{_2}(x\lor y)\thickapprox ( J_{_2}x\land J_{_2}y)\vee( J_{_2}x\land J_{_2}\neg y)\vee ( J_{_2}\neg x\land J_{_2}y) $; \label{BCA:18}
		\item $\Jzero x \thickapprox \Jzero y \;\&\; \Jdue x \thickapprox \Jdue y \;\Rightarrow\; x \thickapprox y$.\label{BCA:quasi}
	\end{enumerate}
\end{definition}

\vspace{5pt}

Bochvar algebras form a proper quasivariety (\cite{FinnGrigolia80}), which will be denoted by $\class{BCA}$, and it is generated by $\textbf{WK}^e$. An immediate example of Bochvar algebras is obtained by any Boolean algebra $\B$ expanded with a new operator $\Jdue=id_B$; in this context $\Jzero$ coincides with the Boolean negation and $\Juno$ is the constant 0.

Axioms (1-8) tells us that (the $\Jdue$-free reduct of) Bochvar algebras have the structure of involutive bisemilattices, as expected since they are generated by $\WKt$, which has $\WK$ as reduct. We can further restrict the class of such involutive bisemilattices, noticing that the quasiequation (FP) $x\approx\neg x$ $\&$ $y\approx\neg y\Rightarrow x\approx y$ holds in $\WK$, therefore it transfers to $\mathbb{Q}(\WKt)$. The quasiequation states that there can be at most one single point of a Bochvar algebra which is a fixpoint for negation. The subquasivariety of $\IBSL$ characterized by (FP) is $\SIBSL$, single-fixpoint involutive bisemilattices (\cite{Bonzio16SL}).

\begin{proposition}[\cite{SMikBochvar}, Proposition 3.5]\label{prop_bca_isbl}
	The $\Jdue$-free reduct of a Bochvar algebra is a $\class{SIBSL}$.
\end{proposition}

As an involutive bisemilattice, the $\Jdue$-free reduct of a Bochvar algebras can be represented as a P\l onka sum of Boolean algebras, guiding us towards a better understanding of the structure of these algebras.

\begin{proposition}[\cite{SMikBochvar}, Lemma 3.9, Remark 3.16]\label{prop.: fatti_BCA}
	Let $\A\in\BCA$ be s.t. $\PL(\mathbb{A})$ is the P\l onka sum decomposition of its $\IBSL$-reduct, with $\mathbb{A} = \pair{\{\A_i\}_{i\in I}, \pair{I,\leq}, \{p_{ij}|i,j\in I, i\leq j\}}$. For $i<j$, all the homomorphisms $p_{ij}$ are surjective and not injective.
\end{proposition}

Notice though that in general there is no way to represent a Bochvar algebra (in the full type) as a P\l onka sum (unless its underlying semilattice direct system is trivial, hence the single resulting fiber is just a Boolean algebra). Once the external operator $\Jdue$ is added to the language the structure of P\l onka sum is broken, since, as the next theorem will state, the application of $\Jdue$ maps any element to the lowest fiber $\A_0$, which is not how operations should be computed in a P\l onka sum. As an example, for an arbitrary fiber $A_i$ and $a\in A_i$, $\Jdue(a\land\neg a)=0\in A_0$, because by Definition \ref{def.: plonka_sum} the constants always inhabit the lowest fiber.

Despite this negative result, the $\Jdue$ operator has a strong interaction with the P\l onka sum representation of the $\IBSL$-reduct of a Bochvar. The application of $\Jdue$ to a fiber of the P\l onka sum representation of a Bochvar algebra produces an isomorphic copy of that fiber inside the bottom fiber. 

\begin{theorem}[\cite{SMikBochvar}, Lemmas 3.9-3.14, Theorems 3.12-3.15], \label{th.: BCA_representation}
	Let $\A\in\BCA$ be s.t. $\PL(\mathbb{A})$ is the P\l onka sum decomposition of its $\IBSL$-reduct, with $\mathbb{A} = \pair{\{\A_i\}_{i\in I}, \pair{I,\leq}, \{p_{ij}|i,j\in I, i\leq j\}}$. The following hold:
	\begin{itemize}
		\item[(i)] for every $i\in I$ and $a\in A_i$, $\Jdue(a)\in p_{0i}^{-1}(a)$; in particular for $a\in A_0, \Jdue(a)=a$;
		
		\item[(ii)] for every $i\in I$, $\Jdue:\A_i\to\mathbf{[0, \Jdue 1_i\mathbf]}$ is a (Boolean) isomorphism (with inverse the restriction of $p_{0i}$ to $\Jdue[\A_i]$), with $\mathbf{[0, \Jdue 1_i\mathbf]}:=\pair{[0, \Jdue 1_i\mathbf],\lor,\land,^*,0,\Jdue 1_i}$ the interval Boolean algebra, where negation $^*$ is defined as $a^*:=\neg a\land\Jdue 1_i$;
		
		\item[(iii)] for every $i \neq j$, $\Jdue 1_i \neq \Jdue 1_j$; in particular, if $i < j$ then $\Jdue 1_j < \Jdue 1_i$.
	\end{itemize}
\end{theorem}

The above results state that the bottom fiber of the $\IBSL$-reduct of a Bochvar algebra contains all the information of the upper fibers, in the form of isomorphic copies of those fibers within the bottom one. The external operator $\Jdue$ maps a fiber $A_i$ into the isomorphic interval Boolean algebra $\mathbf{[0, \Jdue 1_i\mathbf]}$ contained in the bottom fiber $A_0$, and the P\l onka homomorphisms of the underlying semilattice direct system isomorphically sends the interval algebra back to the original fiber $A_i$. 

Furthermore $\Jdue$ keeps track of the position of the fibers within the semilattice direct system, since given two different but comparable indexes $i<j$, the top element of the interval algebra $\Jdue[\A_j]$ will be strictly below the top of $\Jdue[\A_i]$ (notice that while in a Bochvar algebra there are two, possibly different, semilattice orderings to consider, it makes sense to refer to a single ordering within $A_0$ since this is a Boolean algebra).

Finally, we explain in what sense $\BCA$ is the algebraic counterpart of weak Kleene logics. In order to do that, we provide some basics of abstract algebraic logic.

\begin{definition}\label{def.: algebrizzabilità}[\cite{Font16_handbook}, Definition 3.11, Proposition 3.12, Corollary 3.18]
	A finitary logic $\logica{\Lambda}$ is \textit{algebraizable} w.r.t. a class of algebras $\class{K}$ if there exist two maps (called transformers) $\tau\colon Fm\to \mathscr{P}(Eq)$ and $\rho\colon Eq\to \mathscr{P}(Fm)$, respectively from formulae to sets of equations and from equations to sets of formulae, both commuting with substitutions, such that:
	\begin{align*}
		& \text{(ALG1)}\;\;\; \Gamma\vdash_{\logica{\Lambda}}\varphi\iff\tau[\Gamma]\vDash_{\mathsf{K}}\tau(\varphi),\\
		& \text{(ALG4)}\;\;\; \varphi\thickapprox\psi\Dashv\vDash_{\class{K}}\tau\circ \rho(\varphi\thickapprox\psi).
	\end{align*}
	
	\n The above conditions are equivalent to:
	\begin{align*}
		& \text{(ALG2)}\;\;\; \Theta\vDash_{\class{K}}\epsilon\approx\delta\iff\rho[\Theta]\vdash_{\logica{\Lambda}}\rho(\epsilon\approx\delta),\\
		& \text{(ALG3)}\;\;\; \varphi\dashv\vdash_{\logica{\Lambda}}\rho\circ \tau(\varphi).
	\end{align*}
	
	\n If a logic $\logica{\Lambda}$ is algebraizable w.r.t. a class $\class{K}, \mathbb{Q}(\class{K})$ is called its \textit{equivalent algebraic semantics} and it is the largest quasivariety w.r.t. which $\logica{\Lambda}$ is algebraizable.
\end{definition}

While weak Kleene logics have very poor algebraic properties\footnote{Neither $\B$ nor $\PWK$ are even protoalgebraic (see \cite[Theorems 6.1.14-5.5.1]{Bonziobook}), which is the most general class within the so-called Leibniz hierarchy, a classification of logics according to the strength of their relation with their algebraic counterpart.}, their external counterparts are instead algebraizable.

\begin{theorem}[\cite{Ignoranzasevera}, Theorem 35]
	$\B_\mathrm{e}$ is algebraizable w.r.t. $\class{BCA}$ with transformers $\tau(\varphi) := \lbrace\varphi\approx 1\rbrace$ and $\rho(\epsilon\approx\delta) := \lbrace\epsilon\equiv\delta\rbrace$.
\end{theorem}
\begin{theorem}[\cite{bonzio_zamperlin25}, Theorem 7.1]\label{th: algebrizzabilita_PWKe}
	$\PWKe$ is algebraizable w.r.t. $\class{BCA}$ with transformers $\tau(\varphi) := \lbrace\neg\Jzero\varphi\approx 1\rbrace$ and $\rho(\epsilon\approx\delta) := \lbrace\epsilon\equiv\delta\rbrace$.
\end{theorem}

\section{Modal Bochvar and Halld\'{e}n algebras}\label{sec: modal_bochvar_algebras}\label{sct:BocHallalgebras}

In this section we investigate the algebraic counterparts of the modal external weak Kleene logics. We will make of the structure theory of Bochvar algebras. In particular, recall that any Bochvar algebra has a $\SIBSL$-reduct, which is a P\l onka sum of Boolean algebras, and that its underlying semilattice direct system has surjective and non-injective homomorphisms. Furthermore the lowest fiber of the semilattice contains isomorphic copies of all the fibers through the mapping induced by the external operator $\Jdue$.

The modal logics we are going to examine are the \textit{global} counterparts of modal external Bochvar ($\MBeKglobal$, Definition \ref{def: conseguenza logica Bochvar}) and modal external $\PWK$ ($\MPWKglobal$, Definition \ref{def: conseguenza logica pwk}), which are obtained by switching the rule of necessitation with the stronger global necessitation:
$$(\text{GN})\hspace{0.2cm} \phi\vdash\Box\phi$$

\noindent Like in standard modal logic, these global versions are extensions of their local counterparts, and the global logic coincide with the local on theorems, as it is immediate to see. The motivation behind this move will be clear in a moment, first we need to extend the completeness theorems proved for the local counterparts to the new global modal logics.

\begin{theorem}[Completeness]\label{th.: completezza_logiche_globali}
	
	\begin{enumerate}
		\item $\Gamma\vdash_{\MBeKglobal}\varphi$ if and only if $ \Gamma\vDash^{g}_{\MBeK}\varphi$;
		\item $\Gamma\vdash_{\MPWKglobal}\varphi$ if and only if $ \Gamma\vDash^{g}_{\MPWK}\varphi$.
	\end{enumerate}
\end{theorem}
\begin{proof}
	The proof is an adaptation of the completeness theorem for the minimal global modal logic $\logica{K}^g$. We will only show (1); the adaptation for (2) will be clear. \\
	The proof relies on the modification of the canonical model used in Theorem \ref{th: completezza Bochvar modale} (see \cite{bonzio_zamperlin25}). Given a certain $\MBeKglobal$-theory $\Gamma$, change Definition \ref{def: modello canonico} of canonical model $\M=\pair{\mathcal{W},\mathcal{R},v}$ by putting as set of worlds all \textit{local} theories which extend $\Gamma$, that is
	$$\mathcal{W}=\{\Delta\supseteq\Gamma\text{ | }\Delta \text{ is a maximally }\MBeK\text{-consistent theory}\}.$$
	We prove that the Truth Lemma \ref{lemma: truth lemma} holds also for $\M$. Observe that the truth lemma requires the Existence Lemma \ref{lemma: di esistenza} to work, but that lemma can be easily adapted to the case of $\MBeKglobal$, in particular notice that since we are working with \textit{local} theories as points of our canonical model, we can employ the deduction theorem for $\MBeK$ (\cite[Theorem 3.4]{bonzio_zamperlin25}), which otherwise wouldn't hold for its global counterpart. In the case of the global logic though we need to prove that the set $s$ of Lemma \ref{lemma: di esistenza} (as proved in \cite[Lemma 3.20]{bonzio_zamperlin25}) - the set that extends $s^-$ and witnesses the instantiation of the argument of a diamond formula - belongs to the new canonical model, that is $\Gamma\subseteq s$. Recall that $s$ is obtained from $s^-= \{\Jdue\psi\text{ | }\Box\Jdue\psi\in w\}$, with $w\in\mathcal{W}$. By definition of the canonical model, $\Gamma\subseteq w$; since $\Gamma$ is a $\MBeKglobal$-theory, for all $\gamma\in\Gamma, \Jdue\gamma\in\Gamma$ (by $\phi\vdash_{\Be}\Jdue\phi$) and $\Box\Jdue\gamma\in w$ (by (GN)). Therefore, for all $\gamma\in\Gamma,\Box\Jdue\gamma\in w$, hence $\Jdue\gamma\in s^-$. Now we extend $s^-$ to a maximally $\MBeK$-consistent set by Lindenbaum's Lemma \ref{lemma: Lindenbaum}, which we can do, since we are extending it to a theory of the local logic.
	
	Now that we have all the needed lemmas, completeness follows. Let $\Gamma\nvdash_{\MBeKglobal}\phi$, expand it to the least $\MBeKglobal$-theory $\Gamma'\supseteq\Gamma$, for which still $\phi\notin\Gamma'$. Since $\MBeKglobal$ extends $\MBeK$, $\Gamma'$ is also a $\MBeK$-theory, which we can further expand by Lindenbaum's model \ref{lemma: Lindenbaum} to a maximally consistent $\MBeK$-theory $\Delta$. Consider the canonical model $\M=\pair{\mathcal{W},\mathcal{R},v}$ for $\Gamma'$: $\Delta\in\mathcal{W}$, furthermore for all $\Sigma\in\mathcal{W}, \Gamma\subseteq\Gamma'\subseteq\Sigma$, so $\M\vDash_{\MBeK}\gamma$, for all $\gamma\in\Gamma$, but $\phi\notin\Delta$, hence $\M\nvDash_{\MBeK}\phi$. We conclude that $\Gamma\nvDash_{\MBeK}^g\phi$.
	
	For claim (2) it is straightforward to adapt the above strategy to $\MPWKglobal$.
\end{proof}

\n With the completeness results it is now legitimate to employ Kripke semantics to check validity. 

The decision to move to global modal logics is due to their better algebraic behaviour compared to their local versions, which tend to have a weak connection with any algebraic counterpart\footnote{It is a well-known result that the local normal modal logics based on $\logica{S5}$ and weaker systems are not algebraizable, while their global versions are not only algebraizable but even implicative (see e.g. \cite{Font16_handbook}, example 3.61).}. On the contrary, the stronger inferential power of global modal logics allow for a meaningful connection between the logic and a class of algebras, in our case it will be in the sense of algebraizability, which will amount to the possibility to translate the logical consequence and the relative equational consequence of to a certain class of algebras into each other and without any loss of information.

To prove that $\MBeKglobal$ and $\MPWKglobal$ are algebraizable we will employ a purely syntactic characterization of algebraizability (\cite[Theorem 3.21]{Font16_handbook}), which is equivalent to the standard Definition \ref{def.: algebrizzabilità}. This allow us to check the algebraizability of expansions of algebraizable logics directly through their proof theory, before knowing their equivalent algebraic semantics:
\begin{proposition}[\cite{Font16_handbook}, proposition 3.31]\label{prop.: caratterizzazione_sintattica_algebrizzabilita}
	Let $\logica{\Lambda}$ be an algebraizable logic and $\Delta(x,y)$ one of its sets of equivalence formulae. An expansion $\logica{\Lambda}'$ of $\logica{\Lambda}$ that satisfies for each $n$-ary symbol $\lambda$ in the expanded type of $\logica{\Lambda}'$:
	$$\bigcup\limits_{1\leq i\leq n}\Delta(x_i,y_i)\vdash_{\logica{\Lambda}'}\Delta(\lambda x_1, ... x_n,\lambda y_1, ... y_n)$$
	\noindent is algebraizable, and its algebraizability is witnessed by the same transformers of $\logica{\Lambda}$. 
\end{proposition}

\n Since $\MBeKglobal$ and $\MPWKglobal$ are expansions of algebraizable logics, we use the above proposition to prove their algebraizability.
\begin{theorem}
	$\MBeKglobal$ and $\MPWKglobal$ are algebraizable with the same transformers of, respectively, $\Be$ and $\PWKe$.
\end{theorem}
\begin{proof}
	Consider $\MBeKglobal$. By Proposition \ref{prop.: caratterizzazione_sintattica_algebrizzabilita}, to prove our claim we simply need to check whether $x\equiv y\vdash_{\MBeKglobal}\Box x\equiv\Box y$. By completeness Theorem \ref{th.: completezza_logiche_globali}, we can prove semantically that $x\equiv y\vDash^g_{\MBeK}\Box x\equiv\Box y$. Take an arbitrary Bochvar-Kripke model $\M$ s.t. for all $x\equiv y$ is valid. For arbitrary $w\in W$, we have three cases:
	\begin{itemize}
		\item[(i)] If $v(w,\Box x)=\ant$ then $v(w,x)=\ant$, which by $v(w,x\equiv y)=1$ yields $v(w,y)=\ant$, so also $v(w,\Box y)=\ant$. We conclude $v(w,\Box x\equiv\Box y)=1$.
		
		\item[(ii)] If $v(w,\Box x)=1$, $v(w,x)\neq\ant$ and then $v(w,y)\neq\ant$ as well, so $v(w,\Box y)\neq\ant$. Since $v(w,\Box x)=1$, in all successors $w'\in R[w],v(w',x)=1$, moreover we assumed $v(w',x\equiv y)=1$, therefore $v(w',y)=1$, which implies $v(w,\Box y)=1$. We conclude $v(w,\Box x\equiv\Box y)=1$.
		
		\item[(iii)] $v(w,\Box x)=0$ implies $v(w,\Box y)\neq\ant$, as in the previous point. Since $v(w,\Box x)=0$, there is a successors $w'\in R[w]$ s.t. $v(w',x)\neq 1$, hence by the assumption $v(w',x\equiv y)=1$ we also have $v(w',y)\neq 1$, which gives $v(w,\Box y)=0$. We conclude $v(w,\Box x\equiv\Box y)=1$.
	\end{itemize}
	
	\n $\equiv$ is symmetric, so we don't need to prove anything for $y$ on the right-hand side. The case for $\MPWKglobal$ is completely similar, employing Halld\'{e}n-Kripke models instead.
\end{proof}

\n In the above proof the global reading of the logical consequence played a key role, this motivates our shift towards global modal logics.

From a finitary algebraizable logic with finitary transformers it is possible to extract a quasiequational theory for its equivalent algebraic semantics using the algorithm provided by \cite[proposition 3.44]{Font16_handbook}. Since our global modal logics inherited the finitary transformers of their non-modal counterparts, we can apply the mentioned algorithm and obtain the following quasivarieties.

\begin{definition}\label{modal_Bochvar_algebra}
	The class of \textit{modal Bochvar algebras} $\MBCABoch$, in the language $\pair{\lor,\neg,\Jdue,\Box,0,1}$ of type $\langle 2,1,1,1,0,0\rangle$, is the quasivariety axiomatized by the following:
	\begin{itemize}
		\item all the equations and quasiequations of the basis of $\BCA$;
		
		\item $\square(\Jdue\varphi\to\Jdue\psi)\to(\square\Jdue\varphi\to\square\Jdue\psi)\approx 1$;
		
		\item $\neg\Juno\varphi\leftrightarrow\neg\Juno\square \varphi\approx 1$;
		
		\item $\Jdue\square\varphi\to\square\Jdue\varphi\approx 1$;
		
		\item $\Jzero\square\varphi\to\neg\square\Jzero\neg\varphi\approx 1$;
		
		%\item $\Jdue(\Jdue\phi\leftrightarrow\Jdue\psi)\approx 1$ $\&$ $\Jdue(\Jzero\phi\leftrightarrow\Jzero\psi)\approx 1\Rightarrow\Jdue(\phi\equiv\psi)\approx 1$;
		
		%\item $\Jdue\phi\approx 1\Rightarrow\Jdue(\phi\leftrightarrow 1)\approx 1$;
		
		%\item $\Jdue(\phi\leftrightarrow 1)\approx 1\Rightarrow\Jdue\phi\approx 1$;
		
		\item $\phi\approx 1\Rightarrow\Box\phi\approx 1$;
		
		\item $\phi\equiv\psi\approx 1\Rightarrow\phi\approx\psi$.
	\end{itemize}
	
	\n Using the structural properties of $\BCA$, this axiomatization can be equivalently rewritten as:
	
	\begin{itemize}
		\item[(B$_\theta$)] all the equations and quasiequations of the basis of $\BCA$;
		
		\item[(M1)]  $\square(J_2\varphi\to J_2\psi)\preccurlyeq(\square J_2\varphi\to\square J_2\psi)$;
		
		\item[(M2)]  $\neg\Juno\varphi\approx \neg\Juno\square \varphi$;
		
		\item[(MB3)] $\Jdue\square\varphi\preccurlyeq\square\Jdue \varphi$;
		
		\item[(MB4)] $\Jzero\square\varphi\preccurlyeq\neg\square\Jzero\neg\varphi$;
		
		\item[($\tau$-GN)] $\phi\approx 1\Rightarrow\Box\phi\approx 1$.
	\end{itemize}
	
	\n $\MBCABoch$ is the equivalent algebraic semantics of $\MBeKglobal$.
\end{definition}

\n We have made use of the symbol $\preccurlyeq$ despite the structures we are considering have a bisemilattice reduct (fact guaranteed by the presence of all the quasiequational theory of $\BCA$ in the above axiomatization). This notation is not ambiguous though, since each fiber of the P\l onka sum representation of said bisemilattice is a Boolean algebra and, as recalled in Proposition \ref{prop.: fatti_BCA}, the $J_i$ operators map every point in the algebra to the lowest fiber. Therefore as long as we work with external formulae $\alpha,\beta$, an equation $\alpha\to\beta\approx 1$ can be rewritten as $\alpha\preccurlyeq\beta$ by standard properties of Boolean algebras, and similarly $\alpha\leftrightarrow\beta\approx 1$ becomes $\alpha\approx\beta$. We still need to prove that the $\Box$ operator is computed inside its origin fiber, otherwise the above remarks would be illegitimate. We are going to do that after the next definition, to which the previous comments apply.

\begin{definition}\label{modal_Hallden_algebra}
	The class of \textit{modal  Halld\'{e}n algebras} $\MBCAHalld$, in the language $\pair{\lor,\neg,\Jdue,\Box,0,1}$ of type $\langle 2,1,1,1,0,0\rangle$, is the quasivariety axiomatized by the following:
	\begin{itemize}
		\item all the equations and quasiequations of the basis of $\BCA$;
		
		\item $\neg\Jzero(\square(\Jdue\varphi\to\Jdue\psi)\to(\square\Jdue\varphi\to\square\Jdue\psi))\approx 1$;
		
		\item $\neg\Jzero(\neg\Juno\varphi\leftrightarrow \neg\Juno\square \varphi)\approx 1$;
		
		\item $\neg\Jzero(\Box\phi\leftrightarrow\Box\neg\Jzero\phi)\approx 1$;
		
		%\item $\neg\Jzero(\Jdue\phi\leftrightarrow\Jdue\psi)\approx 1$ $\&$ $\neg\Jzero(\Jzero\phi\leftrightarrow\Jzero\psi)\approx 1\Rightarrow\neg\Jzero(\phi\equiv\psi)\approx 1$;
		
		%\item $\neg\Jzero\phi\approx 1\Rightarrow\neg\Jzero(\neg\Jzero\phi\leftrightarrow 1)\approx 1$;
		
		%\item $\neg\Jzero(\neg\Jzero\phi\leftrightarrow 1)\approx 1\Rightarrow\neg\Jzero\phi\approx 1$;
		
		\item $\neg\Jzero\phi\approx 1\Rightarrow\neg\Jzero\Box\phi\approx 1$;
		
		\item $\Jdue(\phi\equiv\psi)\approx 1\Rightarrow\phi\approx\psi$.
	\end{itemize}
	
	\n Since for $\alpha$ external, $\BCA\vDash\Jzero\alpha\approx\neg\alpha$ and deleting what already derives from the quasiequational theory of $\BCA$, this axiomatization can be equivalently rewritten as:
	
	\begin{itemize}
		\item[(B$_\theta$)] all the equations and quasiequations of the basis of $\BCA$;
		
		\item[(M1)]  $\square(J_2\varphi\to J_2\psi)\preccurlyeq(\square J_2\varphi\to\square J_2\psi)$;
		
		\item[(M2)]  $\neg\Juno\varphi\approx \neg\Juno\square \varphi$;
		
		\item[(MP3)] $\neg\Jzero(\Box\phi\leftrightarrow\Box\neg\Jzero\phi)\approx 1$;
		
		\item[($\tau$-GN$'$)] $\neg\Jzero\phi\approx 1\Rightarrow\neg\Jzero\Box\phi\approx 1$.
	\end{itemize}
	
	\n $\MBCAHalld$ is the equivalent algebraic semantics of $\MPWKglobal$.
\end{definition}

\begin{lemma}\label{lem.: MBCA_box_in_fibra}
	Let $\A\in\MBCABoch$ or $\A\in\MBCAHalld$. Consider the P\l onka sum representation of the $\IBSL$-reduct of $\A$ and let $\A_i$ be a fiber of the representation. If $a\in A_i$ then $\Box a\in A_i$.
\end{lemma}
\begin{proof}
	Assume either $\A\in\MBCABoch$ or $\A\in\MBCAHalld$, and $a\in A_i$. Suppose $\Box a\in A_k$. Since both classes share the equation (M2), we have $\A\vDash\neg\Juno x\approx\neg\Juno\Box x$, hence $\neg\Juno a=\neg\Juno\Box a$, moreover by $\IBSL\vdash x\approx\neg\neg x$ we obtain $\Juno a=\Juno\Box a$. By equations of $\BCA$ (which can all easily be checked on the algebra $\WKt$ that generates the quasivariety) we have $\Juno a:=\neg(\Jdue a\lor\Jdue\neg a)=\neg\Jdue(a\lor\neg a)=\neg\Jdue 1_i$. At the same time $\Juno\Box a=\neg(\Jdue\Box a\lor\Jdue\neg\Box a)=\neg\Jdue(\Box a\lor\neg\Box a)=\neg\Jdue 1_k$, where we used again the just mentioned lemma to justify the second to last identity. So $\neg\Jdue 1_i=\neg\Jdue 1_k$, therefore $\Jdue 1_i=\Jdue 1_k$, which allow us to conclude $i=k$ by Theorem \ref{th.: BCA_representation}.
\end{proof}

\n The lemma states that the modal operator is computed fiber by fiber, therefore it was correct to claim that once an element is mapped to the lowest fiber by a $J_i$ it remains there even when prefixed by $\Box$. 

Observe that given any Bochvar algebra the axioms of both $\MBCABoch$ and $\MBCAHalld$ are satisfied when $\Box$ is defined as the identity map, therefore the $\Box$-free reduct of these two classes is precisely $\BCA$.

On its $\{\Jdue,\Box\}$-free reduct a modal Bochvar or Halld\'en algebra is an involutive bisemilattice and as such it can be represented as a P\l onka sum of Boolean algebras. At this point we might hope the if we add box and consider only the $\Jdue$-free reduct this would provide fibers that are modal algebras. Unfortunately this is not the case, but first we remind their definition.

\begin{definition}
	A \textit{Boolean algebra with operators} ($\BAO$), which for simplicity we will call \textit{modal algebra}, is an algebra $\langle \A,\Box\rangle$ s.t. $\A$ is a Boolean algebra and $\Box$ is a unary operation satisfying:
	\begin{itemize}
		\item[(K$_\lor$)] $\Box(x\lor\neg x)\approx x\lor\neg x$;
		\item[(K$_\land$)] $\Box(x\land y)\approx\Box x\land\Box y$.
	\end{itemize}
\end{definition}

\begin{lemma}\label{lem.: K_validi_per_BCA_modale}
	(K$_\land$) is valid in $\MBCABoch$ but fails in $\MBCAHalld$; (K$_\lor$) is valid in $\MBCAHalld$ but fails in $\MBCABoch$.
\end{lemma}
\begin{proof}
	We prove all claims by checking the corresponding translations of the equations on Kripke frames, using algebraizability. By (ALG2), $\MBCABoch\vDash$ (K$_\land$) iff $\vdash_{\MBeKglobal}\Box(x\land y)\equiv(\Box x\land\Box y)$ iff, by completeness Theorem \ref{th.: completezza_logiche_globali}, $\vDash_{\MBeK}^g \Box(x\land y)\equiv\Box x\land\Box y$; it is now a straightforward exercise on Bochvar-Kripke models to prove the last statement. Similarly, $\MBCAHalld\vDash$ (K$_\lor$) iff $\vdash_{\MPWKglobal}\Box(x\lor\neg x)\equiv(x\lor\neg x)$ iff $\vDash_{\MPWK}^g \Box(x\lor\neg x)\equiv(x\lor\neg x)$, and we check the last statement on Halld\'en-Kripke models.
	
	For the invalid equations, by the same reasoning above $\MBCABoch\nvDash$ (K$_\lor$) iff $\nvDash_{\MBeK}^g \Box(x\lor\neg x)\equiv(x\lor\neg x)$, which is witnessed by the following Bochvar-Kripke model: $W=\{w_1,w_2\},Rw_1w_2,v(w_1,x)=1$ and $v(w_2,x)=\ant$. Similarly, $\MBCAHalld\nvDash$ (K$_\land$) iff $\nvDash_{\MPWK}^g \Box(x\land y)\equiv\Box x\land\Box y$, witnessed by the following Halld\'en-Kripke model: $W=\{w_1,w_2\},Rw_1w_2,v(w_1,x)=v(w_1,y)=1$, $v(w_2,x)=0$ and $v(w_2,y)=\ant$.
\end{proof}

We can weaken the definition of modal algebra and produce two larger classes.

\begin{definition}
	Let $\mathbf{B}=\pair{\A,\Box}$ be a structure comprising a Boolean algebra and a unary operation $\Box$. $\mathbf{B}$ is a \textit{meet-modal algebra} ($\BAOmeet$) if it satisfies (K$_\land$); it is a \textit{join-modal algebra} ($\BAOjoin$) if it satisfies (K$_\lor$).
\end{definition}

\n Clearly the intersection of $\BAOmeet$ and $\BAOjoin$ returns precisely all the modal algebras.

Using the fact that box is computed fiber by fiber (lemma \ref{lem.: MBCA_box_in_fibra}) and the valid equations of lemma \ref{lem.: K_validi_per_BCA_modale}, we conclude:

\begin{corollary}\label{cor.: MBCA_fibers_semi_modal_algebras}
	For $\A\in\MBCABoch$, each fiber of the P\l onka sum representation of the $\IBSL$-reduct of $\A$  is, considered in the $\Jdue$-free reduct of $\A$, a meet-modal algebra. For $\A\in\MBCAHalld$, each fiber of the P\l onka sum representation of the $\IBSL$-reduct of $\A$  is, considered in the $\Jdue$-free reduct of $\A$, a join-modal algebra.
\end{corollary}

\n The structure of the fibers of a $\MBCAHalld$ is actually richer.

\begin{lemma}\label{lem.: hallden_bca_modal_fibers}
	Let $\A\in\MBCAHalld$, and $\A_i$ a fiber of the P\l onka sum representation of its $\IBSL$-reduct. $\pair{\A_i,\Box}\in\BAO$.
\end{lemma}
\begin{proof}
	We need to prove that for $a,b\in A_i$ it holds $\Box(a\land b)=\Box a\land\Box b$. We can express syntactically the fact that two elements belong to the same fiber, using the property of $\BCA$ that for indexes $j,k$, if $\Jdue 1_j=\Jdue 1_k$ then $j=k$ (Theorem \ref{th.: BCA_representation}). Since $a,b\in A_i$, $a\lor\neg a=b\lor\neg b=1_i$. Hence we need to prove that $\pair{\A_i,\Box}$ satisfies the quasiequation $\Jdue(x\lor\neg x)\approx\Jdue(y\lor\neg y)\Rightarrow\Box(x\land y)\approx\Box x\land\Box y$, that is $\Jdue(x\lor\neg x)\approx\Jdue(y\lor\neg y)\vDash_{\pair{\A_i,\Box}}\Box(x\land y)\approx\Box x\land\Box y$. By (ALG2) this amounts to $\Jdue(x\lor\neg x)\equiv\Jdue(y\lor\neg y)\vdash_{\MPWKglobal}\Box(x\land y)\equiv\Box x\land\Box y$, which by completeness is the same as $\Jdue(x\lor\neg x)\equiv\Jdue(y\lor\neg y)\vDash_{\MPWK}^g\Box(x\land y)\equiv\Box x\land\Box y$.
	
	Let $\M=\pair{W,R,v}$ be a Halld\'en-Kripke model s.t. $\M\vDash\Jdue(x\lor\neg x)\equiv\Jdue(y\lor\neg y)$. Observe that the last statement means that at each world $\Jdue(x\lor\neg x)$ and $\Jdue(y\lor\neg y)$ have the same value: if $v(w,\Jdue(x\lor\neg x))=1$ then $v(w,x\lor\neg x)=1$ and $v(w,x)\neq\ant$, otherwise $v(w,x)=\ant$. Therefore the assumption forces at each world $x$ and $y$ to be either both classical or both non-classical.
	
	Now to prove the main claim, suppose for $w\in W,v(w,\Box(x\land y))=1$, so $v(w,x)\neq\ant\neq v(w,y)$. For all $w'\in R[w],v(w',x\land y)\neq 0$, and since by hypothesis $v(w',\Jdue(x\lor\neg x)\equiv\Jdue(y\lor\neg y))=1$, either $v(w',x)\neq\ant\neq v(w',y)$ or $v(w',x)=v(w',y)=\ant$. In the former case, since $v(w',x\land y)\neq 0$ it must be $v(w',x)=v(w',y)=1$. Hence in both cases $v(w',x)\neq 0\neq v(w',y)$, so $v(w,\Box x)=v(w,\Box y)=1$. We conclude $v(w,\Box x\land\Box y)=1$. The opposite direction that assumes $v(w,\Box x\land\Box y)=1$ holds in general and doesn't make use of the hypothesis.
	
	Suppose $v(w,\Box x\land\Box y)=0$, which holds if $v(w,x)\neq\ant\neq v(w,y)$ and either $v(w,\Box x)=0$ or $v(w,\Box y)=0$. Therefore it exists $w'\in R[w]$ where $v(w',x)=0$ or $v(w',y)=0$. By hypothesis $v(w',\Jdue(x\lor\neg x)\equiv\Jdue(y\lor\neg y))=1$, hence $v(w',x)\neq\ant$ iff $v(w',y)\neq\ant$, therefore $v(w',x\land y)=0$. We conclude $v(w,\Box(x\land y))=0$. As before, the opposite direction doesn't make use of the hypothesis and it holds in general.
	
	Finally, $v(w,\Box(x\land y))=\ant$ iff ($v(w,x)=\ant$ or $v(w,y)=\ant$) iff ($v(w,\Box x)=\ant$ or $v(w,\Box y)=\ant$) iff $v(w,\Box x\land\Box y)=\ant$.
\end{proof}

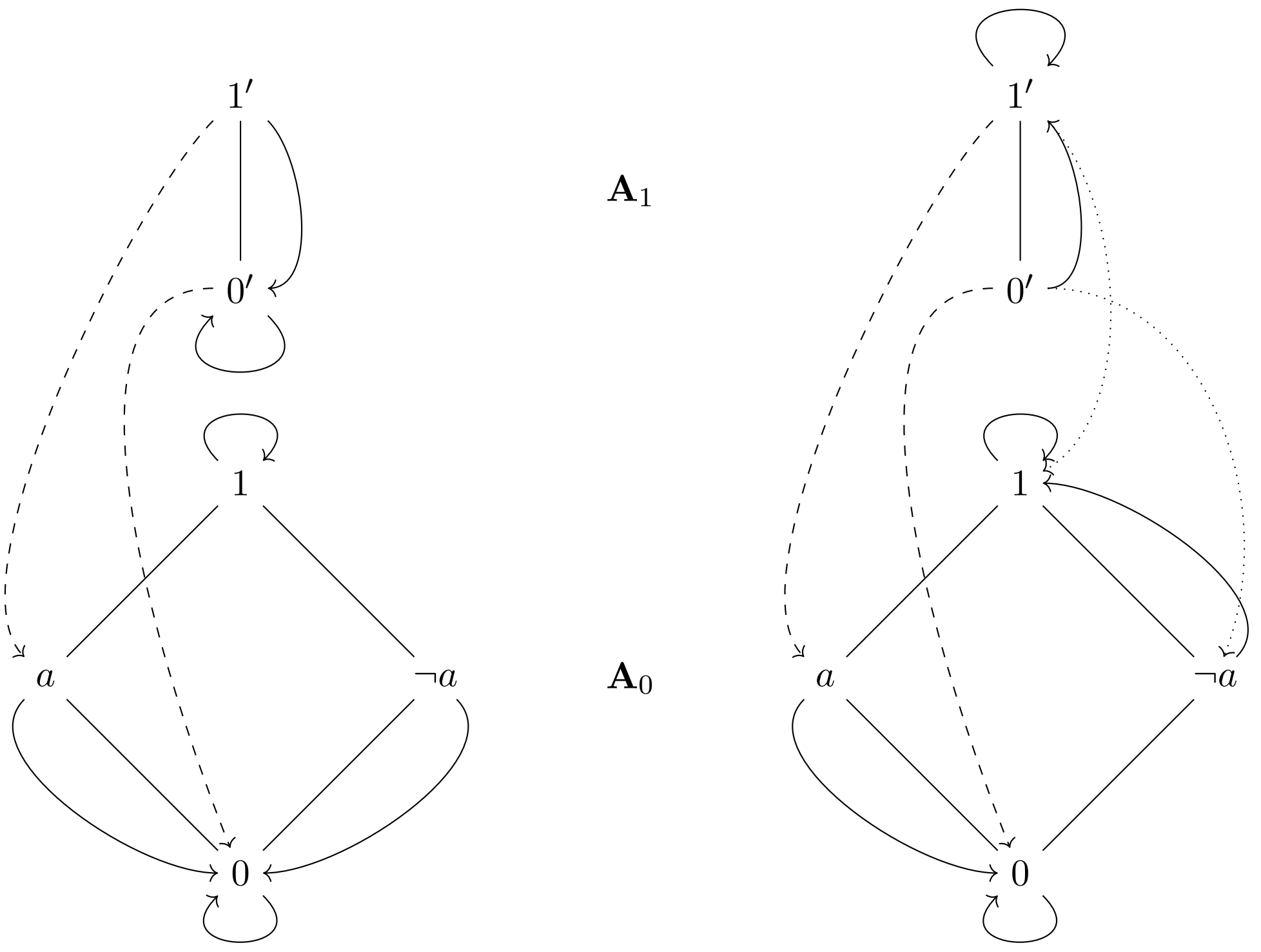
\begin{figure}[h]
	\centering
	\begin{tikzpicture}[scale=1]
		\node (top) at (0,10) {$1'$};
		\node (bot) at (0,8) {$0'$};
		\node (1) at (0,6) {$1$};
		\node (a) at (-2,4) {$a$};
		\node (a') at (2,4) {$\neg a$};
		\node (0) at (0,2) {$0$};
		
		\draw (top) -- (bot);
		\draw (1) -- (a) -- (0) -- (a') -- (1);
		\draw [->] (1) to [out=135,in=45,looseness=5] (1);
		\draw [->] (a) to [out=225,in=180,looseness=0.8] (0);
		\draw [->] (a') to [out=-45,in=0,looseness=0.8] (0);
		\draw [->] (0) to [out=-45,in=-135,looseness=5] (0);
		\draw [->] (top) to [out=-45,in=0,looseness=0.8] (bot);
		\draw [->] (bot) to [out=-45,in=-135,looseness=5] (bot);
		\draw [->,dashed] (top) to [out=225,in=135,looseness=0.5] (a);
		\draw [->,dashed] (bot) to [out=180,in=112,looseness=0.8] (0);
		
		\node (topx) at (8,10) {$1'$};
		\node (botx) at (8,8) {$0'$};
		\node (1x) at (8,6) {$1$};
		\node (ax) at (6,4) {$a$};
		\node (a'x) at (10,4) {$\neg a$};
		\node (0x) at (8,2) {$0$};
		
		\draw (topx) -- (botx);
		\draw (1x) -- (ax) -- (0x) -- (a'x) -- (1x);
		\draw [->] (1x) to [out=135,in=45,looseness=5] (1x);
		\draw [->] (ax) to [out=225,in=180,looseness=0.8] (0x);
		\draw [->] (a'x) to [out=45,in=0,looseness=0.8] (1x);
		\draw [->] (0x) to [out=-45,in=-135,looseness=5] (0x);
		\draw [->] (botx) to [out=0,in=-45,looseness=0.8] (topx);
		\draw [->] (topx) to [out=135,in=45,looseness=5] (topx);
		\draw [->,dashed] (topx) to [out=225,in=135,looseness=0.5] (ax);
		\draw [->,dotted] (topx) to [out=-45,in=28,looseness=0.8] (1x);
		\draw [->,dashed] (botx) to [out=180,in=112,looseness=0.8] (0x);
		\draw [->,dotted] (botx) to [out=0,in=68,looseness=1] (a'x);
		
		\node (topc) at (4,9) {$\A_1$};
		\node (botc) at (4,4) {$\A_0$};
	\end{tikzpicture}
	\caption{A modal Bochvar algebra (left) and a modal Halld\'en algebra (right).}\label{fig.: MBCA}
\end{figure}

The result of corollary \ref{cor.: MBCA_fibers_semi_modal_algebras}  unfortunately does not translate into the possibility to extend the P\l onka homomorphisms to the box operator: considering the P\l onka sum representation of the $\IBSL$-reduct of a modal Bochvar of Halld\'en algebra, we are not guaranteed that the homomorphisms of their semilattice direct system are such w.r.t. $\Box$. The algebras in figure \ref{fig.: MBCA} witness this failure (full arrows represent $\Box$, dashed arrows represent $\Jdue$, which on the bottom fiber is always the identity). In the figure, if we denote by $p_{01}$ the P\l onka (Boolean) homomorphism from the bottom fiber $\A_0$ to the upper one $\A_1$, in the left $\MBCABoch$ observe that $p_{01}(\Box 1)=p_{01}(1)=1'\neq 0'=\Box 1'=\Box p_{01}(1)$. Similarly in the right $\MBCAHalld$ we have $p_{01}(\Box a)=p_{01}(0)=0'\neq 1'=\Box 1'=\Box p_{01}(a)$.

We can actually prove more.
\begin{lemma}\label{lem.: mbca_fallimento_regolarizzazione}
	There is no possible P\l onka sum decomposition of $\MBCABoch$ s.t. its $\Jdue$-free reduct belongs entirely to $\mathcal{R}(\BAOmeet)$. Similarly, There is no possible P\l onka sum decomposition of $\MBCAHalld$ s.t. its $\Jdue$-free reduct belongs entirely to $\mathcal{R}(\BAOjoin)$ nor $\mathcal{R}(\BAO)$.
\end{lemma} 
\begin{proof}
	Let $\MBCABoch^-$ be the $\Jdue$-free reduct of $\MBCABoch$. Suppose by contradiction that $\MBCABoch^-\subseteq\mathcal{R}(\BAOmeet)$. Notice that $\BAOmeet$ is a strongly irregular variety, since by its Boolean base it holds $\BAOmeet\vDash x\land(x\lor y)\approx x$. \cite[Theorem 7.6]{Romanowska86} describes how to obtain from the axiomatization of a strongly irregular variety that of its regularization, in particular $\mathcal{R}(\BAOmeet)\vDash \Box x\land(\Box x\lor y)\approx \Box(x\land(x\lor y))$, so by our assumption also $\MBCABoch\vDash \Box x\land(\Box x\lor y)\approx \Box(x\land(x\lor y))$. By (ALG2) this amounts to $\vdash_{\MBeK}(\Box x\land(\Box x\lor y))\equiv \Box(x\land(x\lor y))$, which by completeness amounts to $\vDash_{\MBeK}(\Box x\land(\Box x\lor y))\equiv \Box(x\land(x\lor y))$. But this fails: consider the Bochvar-Kripke model with $W=\{w_1,w_2\}$ and $w_1Rw_2$ as only relation. Put $v(w_1,x)=v(w_1,y)=v(w_2,x)=1$, while $v(w_2,y)=\ant$; it results $v(w_1,\Box x\land(\Box x\lor y))=1\neq 0=\Box(x\land(x\lor y))$.
	
	The strategy for $\MBCAHalld$ is the same. Suppose by contradiction that the $\Jdue$-free reduct $\MBCAHalld^-\subseteq\mathcal{R}(\BAOjoin)$. This amounts to $\vDash_{\MPWK}(\Box x\land(\Box x\lor y))\approx \Box(x\land(x\lor y))$. This fails: consider the Halld\'en-Kripke model with $W=\{w_1,w_2\}$ and $w_1Rw_2$ as only relation. Put $v(w_1,x)=v(w_1,y)=1,v(w_2,x)=0$, while $v(w_2,y)=\ant$; it results $v(w_1,\Box x\land(\Box x\lor y))=0\neq 1=v(w_1,\Box(x\land(x\lor y)))$.
\end{proof}

\begin{corollary}
	For $\A\in\MBCABoch$ and $\A^-$ its $\Jdue$-free reduct, $\A^-\in\mathcal{R}(\BAOmeet)$ iff $\A^-\vDash\Box x\land(\Box x\lor y)\approx \Box(x\land(x\lor y))$. For $\A\in\MBCAHalld, \A^-\in\mathcal{R}(\BAOjoin)$ iff $\A^-\vDash\Box x\land(\Box x\lor y)\approx \Box(x\land(x\lor y))$. 
\end{corollary}
\begin{proof}
	By the proof of lemma \ref{lem.: mbca_fallimento_regolarizzazione}, we know that if $\A^-\in\mathcal{R}(\BAOmeet)$ then in particular $\A^-\vDash\Box x\land(\Box x\lor y)\approx \Box(x\land(x\lor y))$. For the other direction, let $\A^-\vDash\Box x\land(\Box x\lor y)\approx \Box(x\land(x\lor y))$. Let us abbreviate $x\cdot y:=x\land(x\lor y)$. Observe that an axiomatization of $\BAOmeet$ with at most one irregular equation of the form $\phi(x,y)\approx x$ is given by the equations of $\IBSL$ plus $x\cdot y\approx x$ (which delivers $\class{BA}$), plus the modal equation (K$_\land$) $\Box(x\land y)\approx\Box x\land\Box y$, which is regular. This base meets the criteria to apply \cite[Theorem 7.6]{Romanowska86}, which provides an axiomatization of $\mathcal{R}(\BAOmeet)$ composed by the base of $\IBSL$ and (K$_\land$), plus:
	\begin{itemize}
		\item[(i)] $x\cdot x\approx x$;
		\item[(ii)] $(x\cdot y)\cdot z\approx x\cdot(y\cdot z)$;
		\item[(iii)] $x\cdot(y\cdot z)\approx x\cdot(z\cdot y)$;
		\item[(iv)] $\lambda(x_1,...$ $x_n)\cdot y\approx \lambda(x_1\cdot y,...$ $x_n\cdot y)$;
		\item[(v)] $y\cdot\lambda(x_1,...$ $x_n)\approx y\cdot x_1\cdot ... \cdot x_n$;
	\end{itemize}
	
	\n with $\lambda$ any $n$-ary operation of the type. In $\MBCABoch$ (K$_\land$) holds and this class has a $\IBSL$-reduct, therefore all equations are immediately satisfied except for (iv-v) in the case $\lambda=\Box$. (v) can be checked on Bochvar-Kripke models. (iv) is the critical equation, but we assumed it to be valid in the particular algebra $\A^-$. Therefore $\A^-$ satisfies all the equational base of $\mathcal{R}(\BAOmeet)$, hence it belongs to the class.
	
	The identical reasoning applies for the second claim about $\MBCAHalld$, with the exception that in the equational base for $\mathcal{R}(\BAOjoin)$ we have (K$_\lor$) $\Box(x\lor\neg x)\approx x\lor\neg x$ instead of (K$_\land$) (and both equations in the case of $\mathcal{R}(\BAO)$, the claim about which follows immediately from Lemma \ref{lem.: hallden_bca_modal_fibers} and Theorem \ref{th.: plonka_valid_equations}).
\end{proof}

Despite the fibers of the algebras described don't have the full structure of modal algebras, the bottom fiber is actually a $\BAO$:
\begin{lemma}\label{lem.: bottom_bca_bao}
	Let $\A\in\MBCABoch$ or $\A\in\MBCAHalld$. If $\A_0$ is the bottom fiber of the P\l onka sum representation of the $\IBSL$-reduct of $\A$, $\pair{\A_0,\Box}\in\BAO$.
\end{lemma}
\begin{proof}
	Take $\A\in\MBCABoch$. By corollary \ref{cor.: MBCA_fibers_semi_modal_algebras}, we already know that $\pair{\A_0,\Box}\in\BAOmeet$. To prove that (K$_\lor$) holds, since we are also in a Boolean algebra, $\pair{\A_0,\Box}\vDash x\lor \neg x\approx 1$. By (GN) $\pair{\A_0,\Box}\vDash \Box 1\approx 1$, hence for all $a\in A_0, \Box(a\lor\neg a)=\Box(1)=1=a\lor\neg a$.
	
	Now take $\A\in\MBCAHalld$. Again by corollary \ref{cor.: MBCA_fibers_semi_modal_algebras}, we have $\pair{\A_0,\Box}\in\BAOjoin$. The validity of (K$_\land$) amounts to prove that for $a,b\in A_0,\Box(a\land b)=\Box a\land\Box b$. In general in a $\BCA$ for an element $a$ to belong to the bottom fiber is equivalent to $\Jdue a=a$, therefore our claim can be rephrased as $\Jdue x\approx x,\Jdue y\approx y\vDash_{\MBCAHalld}\Box(x\land y)\approx\Box x\land\Box y$. Since $\MBCAHalld$ is algebraizable, by (ALG2) the previous holds iff $\Jdue x\equiv x,\Jdue y\equiv y\vdash_{\MPWKglobal}\Box(x\land y)\equiv\Box x\land\Box y$ iff, by completeness Theorem \ref{th.: completezza_logiche_globali}, $\Jdue x\equiv x,\Jdue y\equiv y\vDash_{\MPWK}^g\Box(x\land y)\equiv\Box x\land\Box y$. Consider a Halld\'en-Kripke model $\pair{W,R,v}$ s.t. for all $w\in W, v(w,\Jdue x\equiv x)=v(w,\Jdue y\equiv y)=1$ (the case $\ant$ is excluded since we are dealing with external formulae). Take $w\in W$ s.t. $v(w,\Box(x\land y))=1$, then $v(w,x\land y)\neq\ant$ and for all $w'\in R[w],v(w',x\land y)\neq 0$, which holds iff either $v(w',x)=\ant$ or $v(w',y)=\ant$ or $v(w',x)=v(w',y)=1$. The first two cases cannot be realized, since by assumption $v(w',\Jdue x\equiv x)=v(w',\Jdue y\equiv y)=1$, that is $x$ and $y$ receive only classical values at $w'$. Therefore $v(w',x)=v(w',y)=1$, which implies $v(w,\Box x)=v(w,\Box y)=1$, hence $v(w,\Box x\land\Box y)=1$. The opposite direction from the assumption $v(w,\Box x\land\Box y)=1$ is similar and the cases for $0$ are straightforward.    
\end{proof}

Despite the lack of homomorphism that respect $\Box$ between the fibers of the P\l onka sum representation of a $\MBCABoch$, we can extend the isomorphism of Theorem \ref{th.: BCA_representation} to the modal operator. Remember that for a fiber $\A_i$, $1_i:=1^{\A_i}$.

\begin{theorem}\label{th.: isomorfismo_modal_Bochvar_algebra}
	Let $\A\in\MBCABoch$. If $\A_i$ is a fiber of the P\l onka sum representation of the $\IBSL$-reduct of $\A$, $\Jdue$ is an isomorphism of $\BAOmeet$ between $\pair{\A_i,\Box}$ and $\pair{\mathbf{[0,\Jdue 1_i],\Box^\downarrow}}$, where $\Box^\downarrow x:=\Box x\land\Jdue 1_i$.
\end{theorem}
\begin{proof}
	First we prove that that $\pair{\mathbf{[0,\Jdue 1_i]},\Box^\downarrow}\in\BAOmeet$: $\Box^\downarrow(a\land b)=\Box(a\land b)\land\Jdue 1_i=\Box a\land\Box b\land\Jdue 1_i=\Box a\land\Jdue 1_i\land\Box b\land\Jdue 1_i=\Box^\downarrow(a)\land\Box^\downarrow(b)$, where all the identities hold because $\A_0$ is a $\BAO$. 
	
	By Theorem \ref{th.: BCA_representation} we already know that $\Jdue$ is a Boolean isomorphism. For $\Box$ we need to prove $\MBCABoch\vDash\Jdue\Box x\approx\Box^\downarrow\Jdue x$, where $\Box^\downarrow\Jdue x:=\Box\Jdue x\land\Jdue 1_i$ and $1_i$ can be expressed as $x\lor\neg x$, since any value we are going to assign to $x$ belongs to $\A_i$ (we are working inside a Boolean algebra, we can express the local zero in the usual way). Therefore we need to prove $\MBCABoch\vDash\Jdue\Box x\approx\Box\Jdue x\land\Jdue(x\lor\neg x)$, which by algebraizability plus Kripke completeness amounts to $\vDash_{\MBeK}\Jdue\Box x\equiv\Box\Jdue x\land\Jdue(x\lor\neg x)$. Suppose for any Bochvar-Kripke model and $w\in W$, $v(w,\Jdue\Box x)=1$, so $v(w,\Box x)=1$ and $v(w,x)\neq\ant$; the latter yields $v(w,\Jdue(x\lor\neg x))=1$. For all $w'\in R[w],v(w',x)=1$, then $v(w',\Jdue x)=1$, hence $v(w,\Box\Jdue x)=1$. The right-to-left reasoning is similar. Suppose now $v(w,\Jdue\Box x)=0$, then either $v(w,\Box x)=\ant$ or $v(w,\Box x)=0$. In the former case $v(w,x)=\ant$, so $v(w,\Jdue(x\lor\neg x))=0$ and always $v(w,\Box\Jdue x)\neq\ant$, hence $v(w,\Box\Jdue x\land\Jdue(x\lor\neg x))=0$. In the latter, $v(w,x)\neq\ant$ and $v(w,\Jdue(x\lor\neg x))=1$, moreover there is $w'\in R[w]$ s.t. $v(w',x)\neq 1$, then $v(w',\Jdue x)\neq 1$, hence $v(w,\Box\Jdue x)=0$. The right-to-left reasoning is similar. Finally, if $v(w,x)=\ant$ the formula is satisfied, since both of its sides are external formulae.
\end{proof}

\n As in the case for $\BCA$, the above theorem states that bottom fiber of a $\MBCABoch$ stores all the information of the upper fibers, in the form of isomorphic copies. The isomorphism as stated doesn't extend to the case for $\MBCAHalld$, as the right algebra in figure \ref{fig.: MBCA} shows. Nonetheless, Halld\'en algebras still enjoys the property of containing isomorphic copies of the upper fibers, witnessed by a different yet natural map (an example is given in Figure \ref{fig.: MBCA}, where the map is represented by the dotted line).

\begin{theorem}\label{th.: isomorfismo_modal_Hallden_algebra}
	Let $\A\in\MBCAHalld$. If $\A_i$ is a fiber of the P\l onka sum representation of the $\IBSL$-reduct of $\A$, $\neg\Jzero$ is an isomorphism of modal algebras between $\pair{\A_i,\Box}$ and $\pair{\mathbf{[\neg\Jzero 0_i,1]},\Box^\uparrow}$, where $\mathbf{[\neg\Jzero 0_i,1]}$ is the interval Boolean algebra $\pair{[\neg\Jzero 0_i,1],\lor,\neg^\uparrow,\neg\Jzero 0_i,1}$ with negation defined as $\neg^\uparrow x:=\neg x\lor\neg\Jzero 0_i$, while $\Box^\uparrow x:=\Box x\lor\neg\Jzero 0_i$.
\end{theorem}
\begin{proof}
	First we prove that $\pair{\mathbf{[\neg\Jzero 0_i,1]},\Box^\uparrow}\in\BAO$. For $x,y\in[\neg\Jzero 0_i,1]$: 
	\begin{itemize}
		\item[(K$_\lor$)] $\Box^\uparrow(x\lor\neg x)=\Box^\uparrow 1=\Box 1\lor\neg\Jzero 0_i=1\lor\neg\Jzero 0_i=1$, where the first identity holds in virtue of $x\in A_0$, the fourth because $\A_0\in\BAO$;
		\item[(K$_\land$)] $\Box^\uparrow(x\land y)=\Box (x\land y)\lor\neg\Jzero 0_i=(\Box x\land\Box y)\lor\neg\Jzero 0_i=(\Box x\lor\neg\Jzero 0_i)\land(\Box y\lor\neg\Jzero 0_i)=\Box^\uparrow x\land\Box^\uparrow y$, where the second identity holds because $\A_0\in\BAO$.
	\end{itemize}
	
	In order to check that the map is well-defined, first observe that $\neg\Jzero:=\neg\Jdue\neg$ is the composition of two order-reversing maps (the negations) and an order-preserving one ($\Jdue$), so $\neg\Jzero$ is order-preserving as well. Since the images of the bounds of the domain correspond to those of the codomain (for the upper bound, observe that $\neg\Jzero(x\lor\neg x)\approx 1$ is valid in $\BCA$), then $[\neg\Jzero 0_i,1]$ is the correct codomain.
	
	The map $\neg\Jzero$ preserves Boolean operations since it is the composition of a Boolean homomorphism ($\Jdue$) with two dual Boolean homomorphisms ($\neg$). We are left to prove that it preserves box as well, which amounts to the validity of the equation $\neg\Jzero\Box x\approx\Box\neg\Jzero x\lor\neg\Jzero(x\land\neg x)$ in $\MBCAHalld$. Using algebraizability of $\MBCAHalld$, we can equivalently prove that $\neg\Jzero\Box x\equiv\Box\neg\Jzero x\lor\neg\Jzero(x\land\neg x)$ is a theorem of $\MPWKglobal$, hence of $\MPWK$ as well. We prove this semantically using Kripke completeness of Theorem \ref{th: completezza MPWK}. Take an arbitrary Halld\'en-Kripke model and $w\in W$ s.t. $v(w,\neg\Jzero\Box x)=1$, so $v(w,\Jzero\Box x)=0$, which holds when either $v(w,\Box x)=1$ or $v(w,\Box x)=\ant$. In the latter case $v(w,x)=\ant$, hence $v(w,x\land\neg x)=\ant$, so $v(w,\neg\Jzero(x\land\neg x))=1$ and $v(w,\Box\neg\Jzero x\lor\neg\Jzero(x\land\neg x))=1$. Instead if $v(w,\Box x)=1$ then for all $w'\in R[w],v(w',x)\neq 0$, so $v(w',\neg\Jzero x)=1$. Since always $v(w,\neg\Jzero x)\neq\ant$, we obtain $v(w,\Box\neg\Jzero x)=1$, then $v(w,\Box\neg\Jzero x\lor\neg\Jzero(x\land\neg x))=1$.
	
	For the other direction, suppose $v(w,\Box\neg\Jzero x\lor\neg\Jzero(x\land\neg x))=1$, which holds if either $v(w,\Box\neg\Jzero x)=1$ or $v(w,\neg\Jzero(x\land\neg x))=1$. In the former case, we have two subcases: if $v(w,x)=\ant$ then $v(w,\Box x)=\ant$; on the other hand if $v(w,x)\neq\ant$ we obtain $v(w,\Box x)=1$, since for all $w'\in R[w],v(w',\neg\Jzero x)=1$ (by $v(w,\Box\neg\Jzero x)=1$ and using the fact that we are evaluating an external formula), hence $v(w',\Jzero x)=0$, that is $v(w',x)\neq 0$. In both subcases $v(w,\Box x)\neq 0$, therefore $v(w,\neg\Jzero\Box x)=1$. Instead in the case that $v(w,\neg\Jzero(x\land\neg x))=1$, we have $v(w,\Jzero(x\land\neg x))=0$, that is $v(w,x\land\neg x)\neq 0$, which can happen only if $v(w,x)=\ant$, but then $v(w,\Box x)=\ant$, so $v(w,\neg\Jzero\Box x)=1$. Since we are dealing only with external formulae we need no other possibilities to consider and the claim is proved.
\end{proof}

\section{Duality between algebras and frames}\label{sct:representation}

The final result we provide in this work is a duality between modal Bochvar algebras and Kripke frames, in the style of the standard J\'onsson-Tarski duality that states that every Boolean algebras with operators can be embedded in the full complex algebra of its ultrafilter frame (\cite[Theorem 5.43, p.290]{Blue_book}). 

\subsection{Representation of modal Bochvar algebras}

As J\'onsson-Tarski representation theorem employs Boolean ultrafilters to build the dual object for a modal algebra, that is its so-called ultrafilter frame, the first item that we need in our proof is a suitable notion of ultrafilter for $\MBCABoch$.

\begin{definition}
	Let $\A$ be an algebra with a $\BCA$-reduct. $U\subseteq A$ is a \textit{paracomplete external ultrafilter} ($b$-ultrafilter) if $\Jdue[U]$ is a Boolean ultrafilter of the bottom fiber $\A_0$ and, for all $x\in A$, 
	$$x\in U\Leftrightarrow\Jdue x\in U.$$ 
	The set of all $b$-ultrafilters over $\A$ will be denoted by $\textit{Uf}_B(\A)$.
\end{definition}

\noindent In the following, by ultrafilter we intend a standard Boolean ultrafilter. The following lemma shows that $b$-ultrafilters behaves like lattice filters even when outside of the bottom fiber of a modal Bochvar algebra.

\begin{lemma}\label{lem.:jultrafilter_basics}
	For any $\A\in\BCA, U\in \textit{Uf}_B(\A)$ and $x,y\in A$ the following hold:
	\begin{itemize}
		\item if $x,y\in U$ then $x\land y\in U$;
		\item if $x,y\in A_i, x\in U$ and $x\leq y$ then $y\in U$.
	\end{itemize}
\end{lemma}
\begin{proof}
	Let $x,y\in U$, then $\Jdue x,\Jdue y\in\Jdue[U]$, which is an ultrafilter, hence $\Jdue x\land\Jdue y=\Jdue(x\land y)\in\Jdue[U]$. So $x\land y\in U$.
	
	Now suppose that $x,y\in A_i$ and $x\leq y$. By Theorem \ref{th.: BCA_representation}, $\Jdue$ is a Boolean isomorphism, in particular it is order-preserving, therefore $\Jdue x\leq\Jdue y$. If $x\in U$ then $\Jdue x\in \Jdue[U]$, then, since $\Jdue[U]$ is an ultrafilter, $\Jdue y\in \Jdue[U]$ and $y\in U$.
\end{proof}

We now list a few simple but useful facts about $b$-ultrafilters, which are going to be employed abundantly later, in the easiest cases without mention.

\begin{lemma}\label{lem.:jultrafilter_properties}
	Let $\A\in\BCA$ and $\mathbf{I}$ the semilattice of indexes of its P\l onka sum representation. For all $U\in \textit{Uf}_B(\A), i\in I$ and $x,y\in A$ the following hold:
	\begin{itemize}
		\item[(a)] $U\cap A_0\neq\emptyset$;
		\item[(b)] $0_i\notin U$;
		\item[(c)] if $x=\neg x, x\notin U$;
		\item[(d)] if $U\cap A_i\neq\emptyset$ then $1_i\in U$;
		\item[(e)] if $U\cap A_i\neq\emptyset$, for all $x\in A_i$ either $x\in U$ or $\neg x\in U$;
		\item[(f)] if $U\cap A_{i\lor j}\neq\emptyset$ then $U\cap A_{i}\neq\emptyset$;
		\item[(g)] if $U\cap A_{i}\neq\emptyset\neq U\cap A_{j}$ then $U\cap A_{i\lor j}\neq\emptyset$;
		\item[(h)] if $x\in U\cap A_{i}, i\leq j$ and $U\cap A_{j}\neq\emptyset$ then $p_{ij}(x)\in U$.
	\end{itemize}
\end{lemma}
\begin{proof}
	(a) is an immediate consequence of the fact that $\Jdue[U]$ is an ultrafilter over $\A_0$.
	
	For (b), $\Jdue 0_i=\Jdue(x\land\neg x)=\Jdue x\land\Jdue\neg x=0$ for $x\in A_i$, which clearly cannot belong to $\Jdue[U]$, hence to $U$ as well. When $x$ is the fix-point for negation (c) it is just a special case of the previous point.
	
	For (d), let $x\in U$ for some $x\in A_i$. Since $x\leq 1_i$, by the upward closure of Lemma \ref{lem.:jultrafilter_basics} also $1_i\in U$.
	
	For (e), suppose $U\cap A_i\neq\emptyset$, then by the previous point $1_i\in U$. We have then, for $x\in A_i, \Jdue 1_i=\Jdue(x\lor\neg x)=\Jdue x\lor\neg\Jdue x\in\Jdue[U]$, which is an ultrafilter, hence either $\Jdue x\in\Jdue[U]$ or $\Jdue\neg x\in\Jdue[U]$, therefore $x\in U$ or $\neg x\in U$.
	
	For (f), suppose $U\cap A_{i\lor j}\neq\emptyset$, then by point (d) we have $1_{i\lor j}\in U$, so $\Jdue 1_{i\lor j}\in U$. Since $\Jdue 1_{i\lor j}\leq\Jdue 1_{i}$ (as it can be checked on $\WKt$), by upward closure $\Jdue 1_{i}\in U$, hence $1_{i}\in U$.
	
	Now, for (g), the assumptions $U\cap A_{i}\neq\emptyset\neq U\cap A_{j}$ imply $1_i,1_j\in U$, using previous points. Hence $\Jdue 1_i,\Jdue 1_j\in U$, then $\Jdue 1_i\land\Jdue 1_j=\Jdue(1_i\land 1_j)\in U$, therefore $1_i\land 1_j\in U$, which proves the claim since $1_i\land 1_j\in A_{i\lor j}$.
	
	Finally, for (h), assume $x\in U\cap A_{i}, i\leq j$ and $U\cap A_{j}\neq\emptyset$. From $x\in U$ we have $\Jdue x\in U$, while by $U\cap A_{j}\neq\emptyset$ gives $1_j\in U$. By Lemma \ref{lem.:jultrafilter_basics}, the previous facts imply $\Jdue x\land 1_j\in U$. By the definition of P\l onka sum, the term $\Jdue x\land 1_j$ is computed within the fiber $\A_j$ as $p_{0j}(\Jdue x)\land 1_j$, therefore $p_{0j}(\Jdue x)\land 1_j\in U$ implies $p_{0j}(\Jdue x)\in U$, by Lemma \ref{lem.:jultrafilter_basics}. Observe now that $p_{0i}(\Jdue x)=x$, since when appropriately restricted $p_{0i}$ is the inverse of the map $\Jdue$ with domain $A_i$ (Theorem \ref{th.: BCA_representation}). Therefore $p_{ij}(x)=p_{ij}\circ p_{0i}(\Jdue x)=p_{0j}(\Jdue x)\in U$ (because $p_{ij}\circ p_{0i}=p_{0j}$ by definition of P\l onka sum).
\end{proof}

\noindent Notice in particular that the previous lemma states that $b$-ultrafilters contain no local bottom and, if it exists, no fix-point for negation. Moreover Lemma \ref{lem.:jultrafilter_basics} together with points (b) and (e) of Lemma \ref{lem.:jultrafilter_properties} yields the following corollary.

\begin{corollary}\label{cor.:jultrafilters_boolean}
	Let $\A\in\BCA$ and $U\in \textit{Uf}_B(\A)$. If $U\cap A_i\neq\emptyset$ then $U\cap A_i$ is an ultrafilter over $\A_i$.
\end{corollary}

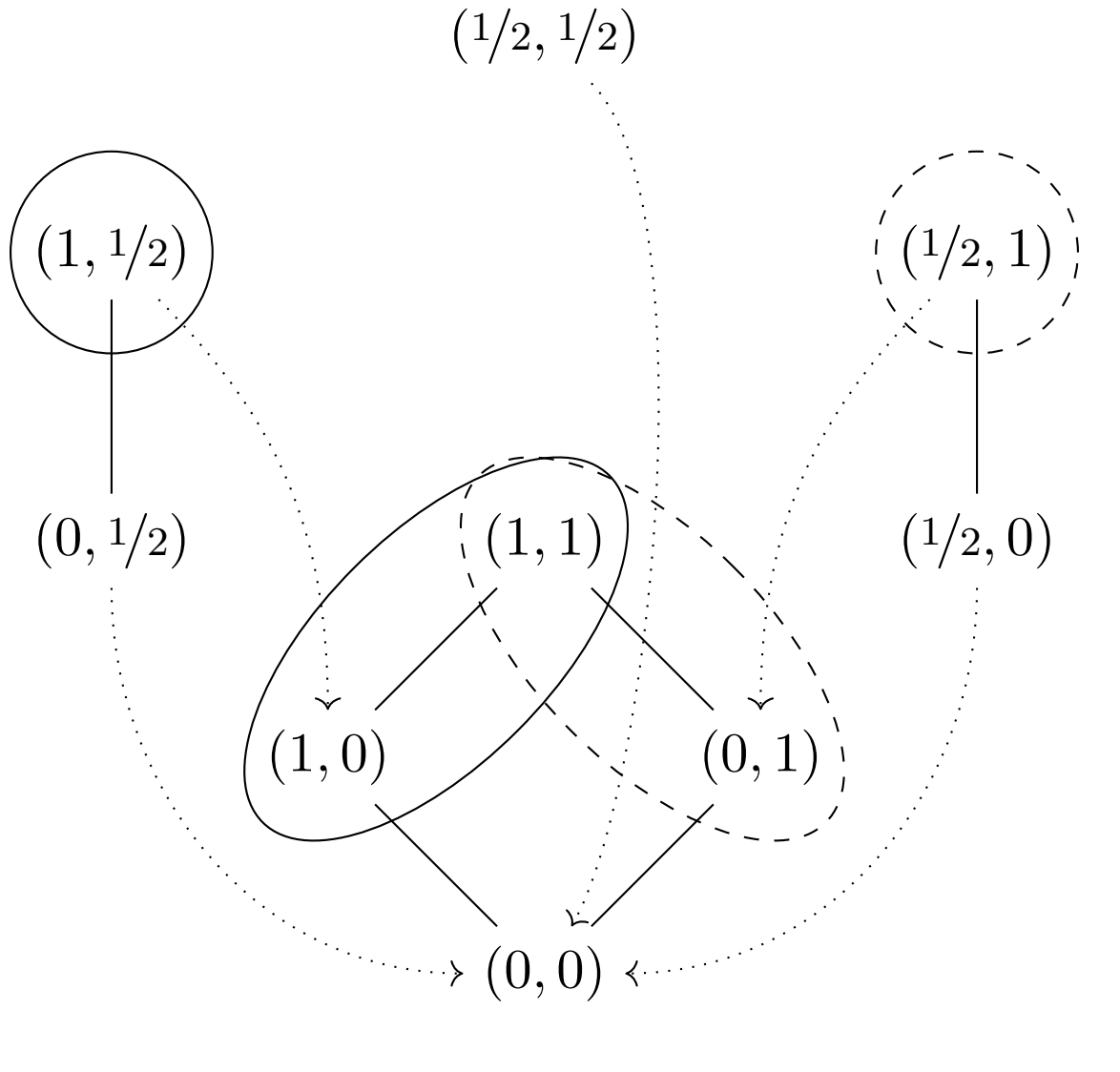
\begin{figure}[h]
	\centering
	\begin{tikzpicture}[scale=1]
		\node (nn) at (0,6.5) {$(\ant,\ant)$};
		\node (1n) at (-3,5) {$(1,\ant)$};
		\node (0n) at (-3,3) {$(0,\ant)$};
		\node (n1) at (3,5) {$(\ant,1)$};
		\node (n0) at (3,3) {$(\ant,0)$};
		\node (11) at (0,3) {$(1,1)$};
		\node (10) at (-1.5,1.5) {$(1,0)$};
		\node (01) at (1.5,1.5) {$(0,1)$};
		\node (00) at (0,0) {$(0,0)$};
		
		\draw (00) -- (01) -- (11) -- (10) -- (00);
		\draw (0n) -- (1n);
		\draw (n0) -- (n1);
		
		\draw (-0.75,2.25) circle [x radius=0.8cm, y radius=1.7cm, rotate=-45];
		\draw[dashed] (0.75,2.25) circle [x radius=0.8cm, y radius=1.7cm, rotate=45];
		\draw (-3,5) circle [radius=0.7cm];
		\draw[dashed] (3,5) circle [radius=0.7cm];
		
		\draw [->,dotted] (1n) to [out=-45,in=90,looseness=1] (10);
		\draw [->,dotted] (0n) to [out=-90,in=180,looseness=1] (00);
		\draw [->,dotted] (n1) to [out=-135,in=90,looseness=1] (01);
		\draw [->,dotted] (n0) to [out=-90,in=0,looseness=1] (00);
		\draw [->,dotted] (nn) to [out=-45,in=60,looseness=0.5] (00);
	\end{tikzpicture}
	\caption{$\WKt^2$ with its two $b$-ultrafilters, represented by circles drawn with two different outlines (dotted lines represent $\Jdue$).}\label{fig.: external_ultrafilter}
\end{figure}

\n From the previous facts it is easy to conclude that for $\A\in\MBCABoch$ there is a bijective correspondence between its external ultrafilters and the Boolean ultrafilters of its bottom fiber $\A_0$. As an example, Figure \ref{fig.: external_ultrafilter} shows the two $b$ -ultrafilters of $\WKt^2$ (represented according to its P\l onka sum decomposition).

Next we need to build a Kripke frame given any modal Bochvar algebra.

\begin{definition}\label{def.:ultrafilter_frame}
	Given $\A\in\MBCABoch$, its $b$-ultraframe is the Kripke frame $\mathfrak{Uf}_B(\A)=\pair{\textit{Uf}_B(\A),Q_A}$, with the accessibility relation defined as
	\begin{align*}
		FQ_AG\Leftrightarrow\text{ for all }x\in A,\Box x\in F\text{ implies }x\in G.
	\end{align*}
\end{definition}

We need to be able to go in the other direction as well, that is from a Kripke to obtain a modal Bochvar algebra.

\begin{definition}\label{def.:bochvar_complex_alg}
	Given a Kripke frame $\mathfrak{F}=\pair{W,R}$, its \textit{full complex Bochvar algebra} is $\mathbf{Cm}_B(\mathfrak{F})=\pair{\{1,\ant,0\}^W,\lor,\neg,\Jdue,m_R,1,0}$, with the operations defined as
	\begin{align*}
		(f\lor g)(w):=f(w)\lor^{\WKt}g(w),\\
		(\neg f)(w):=\neg^{\WKt} f(w),\\
		(\Jdue f)(w):=\Jdue^{\WKt} f(w),\\
		1(w)=1\text{ for all $w\in W$},\\
		0(w)=0\text{ for all $w\in W$},
	\end{align*}
	while
	$$
	(m_R f)(w)=\begin{cases} 
		1 & \text{if }f(w)\neq\ant\text{ and }R[w]\subseteq f^{-1}[\{1\}];\\
		\ant & \text{if }f(w)=\ant;\\
		0 & \text{otherwise.}	
	\end{cases}
	$$
	
	A \textit{complex Bochvar algebra} of $\mathfrak{F}$ is any subalgebra of $\mathbf{Cm}_B(\mathfrak{F})$.	
\end{definition}

\begin{lemma}\label{lem.:complex_alg_is_bochvar}
	For any Kripke frame $\mathfrak{F},\mathbf{Cm}_B(\mathfrak{F})\in\MBCABoch$.
\end{lemma}
\begin{proof}
	Let $\A=\mathbf{Cm}_B(\mathfrak{F})$ for some Kripke frame $\mathfrak{F}=\pair{W,R}$. That the $\Box$-free reduct of $\A$ is a Bochvar algebra is immediate by its definition, since $\A$ is isomorphic to the power $\WKt^{|W|}$.
	
	Next we prove that $\A$ satisfies the modal axioms of $\MBCABoch$. Observe that since all the considered equations involve only external formulae, the cases for $\ant$ will be discarded. Moreover, being external formulae, both sides of the equations are evaluated as objects belonging to the bottom fiber of the P\l onka representation and as such can be compared.
	
	\vspace{0.2cm}
	
	(M1) Suppose $\Box(\Jdue f\to\Jdue g)(w)=\Box(\neg\Jdue f\lor\Jdue g)(w)=1$, which amounts to $R[w]\subseteq(\neg\Jdue f\lor\Jdue g)^{-1}[\{1\}]$ and $\Box(\Jdue f)(w)=1$, so $R[w]\subseteq(\Jdue f)^{-1}[\{1\}]$. For arbitrary $w'\in R[w],(\neg\Jdue f\lor\Jdue g)(w')=1$, which holds iff either $(\neg\Jdue f)(w')=1$ or $(\Jdue g)(w')=1$ iff $f(w')\neq 1$ or $g(w')=1$. At the same time $(\Jdue f)(w')=1$ iff $f(w')=1$, therefore $g(w')=1$ iff $(\Jdue g)(w')=1$. Hence $\Box(\Jdue g)(w)=1$, that is $(\Box\Jdue f\to\Box\Jdue g)(w)=1$.
	
	\vspace{0.2cm}
	
	(M2) $(\neg\Juno f)(w)=1$ iff $(\Juno f)(w)=0$ iff $f(w)\neq\ant$ iff $(\Box f)(w)\neq\ant$ iff $(\Juno\Box f)(w)=0$ iff $(\neg\Juno\Box f)(w)=1$.
	
	\vspace{0.2cm}
	
	(MB3) $(\Jdue\Box f)(w)=1$ iff $(\Box f)(w)=1$ iff $f(w)\neq\ant$ and $R[w]\subseteq f^{-1}[\{1\}]$. Since $f^{-1}[\{1\}]=(\Jdue f)^{-1}[\{1\}]$, $R[w]\subseteq (\Jdue f)^{-1}[\{1\}]$, therefore $(\Box\Jdue f)(w)=1$, because always $(\Jdue f)(w)\neq\ant$.
	
	\vspace{0.2cm}
	
	(MB4) $(\Jzero\Box f)(w)=1$ iff $(\Box f)(w)=0$ iff $f(w)\neq\ant$ and $R[w]\nsubseteq f^{-1}[\{1\}]$. Notice that $f(w)=1$ iff $(\neg f)(w)=0$ iff $(\Jzero\neg f)(w)=1$, therefore $f^{-1}[\{1\}]=(\Jzero\neg f)^{-1}[\{1\}]$. Hence $R[w]\nsubseteq(\Jzero\neg f)^{-1}[\{1\}]$, so $(\Box\Jzero\neg f)(w)=0$ iff $(\neg\Box\Jzero\neg f)(w)=1$.
	
	\vspace{0.2cm}
	
	($\tau$-GN) For all $w\in W,R[w]\subseteq 1^{-1}[\{1\}]=W$, so $\Box 1=1$.
\end{proof}

We are ready to prove our main claim.

\begin{theorem}[Representation]\label{th.:MBCA_representation}
	Let $\A\in\MBCABoch$. There exists an embedding $\iota:\A\to\mathbf{Cm}_B(\mathfrak{Uf}_B(\A))$ given by, for all $x\in A$ and $F\in \textit{Uf}_B(\A)$:
	$$
	\iota(x)(F)=\begin{cases} 
		1 & \text{if }x\in F;\\
		\ant & \text{if }F\cap A_x=\emptyset;\\
		0 & \text{otherwise.}
	\end{cases}
	$$
\end{theorem}
\begin{proof}
	To prove that the map $\iota$, which is clearly well-defined, is injective we consider for $x\neq y\in A$ two cases. If $\Jdue x\neq\Jdue y$ then there exists an ultrafilter $U$ over $\A_0$ which separates the two points, that is, w.l.o.g., $\Jdue x\in U$ and $\Jdue y\notin U$. We can expand $U$ to the $b$-ultrafilter $V=\bigcup_{z\in U} \Jdue^{-1}[z]$, for which $x\in V$ and $y\notin V$. Hence $\iota(x)(V)=1\neq\iota(y)(V)$, so $\iota(x)\neq\iota(y)$. If $\Jdue x=\Jdue y$, it must be the case that $x\in A_i\neq A_j\ni y$, because if $x$ and $y$ belonged to the same fiber their images through $\Jdue$ cannot coincide, since $\Jdue$ restricted to that fiber is an isomorphism. Since $A_i\neq A_j$, also $\Jdue 1_i\neq\Jdue 1_j$ (Theorem \ref{th.: BCA_representation}). Now consider an ultrafilter $U$ s.t. $\Jdue 1_i\in U$ and $\Jdue 1_j\notin U$ and extend it to the $b$-ultrafilter $V$ defined in the previous case. By $\Jdue 1_j\notin U$ we have $1_j\notin V$, therefore $V\cap A_j=\emptyset$, otherwise $V$ should contain at least $1_j$ since it is upward closed by Lemma \ref{lem.:jultrafilter_basics}. We conclude $\iota(x)(V)\neq\ant=\iota(y)(V)$.
	
	Now we show that $\iota$ is a homomorphism.
	
	\vspace{0.2cm}
	
	($\lor$) Suppose $\iota(x\lor y)(F)=1$, which holds iff $x\lor y\in F$ iff $\Jdue(x\lor y)\in F$. By $\Jdue(x\lor y)\leq\Jdue x\lor\Jdue y$ (easy to check on $\WKt$), we have $\Jdue x\lor\Jdue y\in F$, so either $\Jdue x\in F$ or $\Jdue y\in F$, therefore $x\in F$ or $y\in F$. From $F\cap A_{x\lor y}\neq\emptyset$, Lemma \ref{lem.:jultrafilter_properties}, point (f) yields $F\cap A_{x}\neq\emptyset$ and $F\cap A_{y}\neq\emptyset$, then either $\iota(x)(F)=1$ and $\iota(y)(F)\neq\ant$ or $\iota(y)(F)=1$ and $\iota(x)(F)\neq\ant$. Hence $(\iota(x)\lor\iota(y))(F)=1$.
	
	Next, suppose $(\iota(x)\lor\iota(y))(F)=1$, that is either $\iota(x)(F)=1$ and $\iota(y)(F)\neq\ant$, or $\iota(x)(F)\neq\ant$ and $\iota(y)(F)=1$. W.l.o.g., let us assume the first disjunct, which is equivalent to $x\in F$ and $F\cap A_y\neq\emptyset$. By Lemma \ref{lem.:jultrafilter_properties}, point (g) we have $F\cap A_{x\lor y}\neq\emptyset$, which then by point (h) gives $p_{x,x\lor y}(x)\in F$. By Lemma \ref{lem.:jultrafilter_basics}, $p_{x,x\lor y}(x)\lor p_{y,x\lor y}(y)\in F$, but this is precisely the denotation of $x\lor y$ by definition of P\l onka sum. Since we have proved $x\lor y\in F$, this is enough to state $\iota(x\lor y)(F)=1$.
	
	Finally assume $\iota(x\lor y)(F)=\ant$, which holds iff $F\cap A_{x\lor y}=\emptyset$. Using Lemma \ref{lem.:jultrafilter_properties}, points (f-g) by contraposition, the latter holds iff either $F\cap A_x=\emptyset$ or $F\cap A_y=\emptyset$ iff $\iota(x)(F)=\ant$ or $\iota(y)(F)=\ant$ iff $(\iota(x)\lor\iota(y))(F)=\ant$.
	
	\vspace{0.2cm}
	
	($\neg$) Suppose $\iota(\neg x)(F)=1$, which holds iff $\neg x\in F$ iff $F\cap A_{x}\neq\emptyset$ and $x\notin F$ (by Lemma \ref{lem.:jultrafilter_properties}) iff $\iota(x)(F)=0$ iff $(\neg\iota(x))(F)=1$. Next assume $\iota(\neg x)(F)=\ant$, which holds iff $F\cap A_{\neg x}=\emptyset$ iff $F\cap A_{x}=\emptyset$ iff $\iota(x)(F)=\ant$ iff $(\neg\iota(x))(F)=\ant$.
	
	\vspace{0.2cm}
	
	($\Jdue$) Suppose $\iota(\Jdue x)(F)=1$, which holds iff $\Jdue x\in F$. Now since $F$ is a $b$-ultrafilter, $\Jdue x\in F$ iff $x\in F$, hence $F\cap A_x\neq\emptyset$. The last two facts are equivalent to $(\Jdue\iota(x))(F)=1$. Since $(\Jdue\iota(x))(F)$ and $\iota(\Jdue x)(F)$ can never return $\ant$, this is enough to prove the case.
	
	\vspace{0.2cm}
	
	($\Box$) We want to prove $\iota(\Box x)=m_{Q_A}\iota(x)$. Let us first make explicit the clauses of $m_{Q_A}\iota(x)$. For $F\in \textit{Uf}_B(\A)$,
	$$
	(m_{Q_A}\iota(x))(F)=\begin{cases} 
		1 & \text{if }\iota(x)(F)\neq\ant\text{ and }Q_A[F]\subseteq \iota(x)^{-1}[\{1\}];\\
		\ant & \text{if }\iota(x)(F)=\ant;\\
		0 & \text{otherwise}	
	\end{cases}
	$$
	which can be equivalently rephrased as
	$$
	(m_{Q_A}\iota(x))(F)=\begin{cases} 
		1 & \text{if }F\cap A_x\neq\emptyset\text{ and }Q_A[F]\subseteq\{G\in \textit{Uf}_B(\A)|x\in G\};\\
		\ant & \text{if }F\cap A_x=\emptyset;\\
		0 & \text{otherwise.}	
	\end{cases}
	$$
	
	Suppose that $\iota(\Box x)(F)=1$, therefore $\Box x\in F$. By Lemma \ref{lem.: MBCA_box_in_fibra}, we know that every fiber of the P\l onka representation of $\A$ is closed under $\Box$, therefore $A_{\Box x}=A_x$, hence $F\cap A_x\neq\emptyset$. Since $\Box x\in F$, Definition \ref{def.:ultrafilter_frame} implies that for all $G\in Q_A[F],x\in G$, so $Q_A[F]\subseteq\{G\in \textit{Uf}_B(\A)|x\in G\}$. This amounts to $(m_{Q_A}\iota(x))(F)=1$.
	
	Now suppose that $(m_{Q_A}\iota(x))(F)=1$, so $Q_A[F]\subseteq\{G\in \textit{Uf}_B(\A)|x\in G\}$ and $F\cap A_x\neq\emptyset$. Since $A_x=A_{\Box x}$, all we are left to show in order to prove $\iota(\Box x)(F)=1$ is $\Box x\in F$. Assume the contrary, that is $\Box x\notin F$, and consider the sets
	$$
	H=\{\Jdue y\text{ | }\Box\Jdue y\in F\},
	$$$$
	H'=H\cup\{\neg\Jdue x\}.
	$$
	Observe that $H'$ belongs entirely to the bottom fiber $\A_0$ of the P\l onka representation. To demonstrate that the Boolean filter generated by $H'$ is non-trivial, we need to prove that no finite meet of elements of $H'$ is equal to 0. To do that take $z$ finite meet of elements of $H$ and suppose towards contradiction that $z\land\neg\Jdue x=0$. Using basic properties of Boolean algebras this implies $z\leq\Jdue x$, then $\Box z\leq\Box\Jdue x$ by monotonicity of $\Box$ (recall that $\A_0$ is a modal algebra, Lemma \ref{lem.: bottom_bca_bao}). By definition $z$ is of the form $\Jdue y_1\land...\land\Jdue y_n$ for $\Jdue y_1,...\Jdue y_n\in H$, therefore $\Box\Jdue y_1,...\Box\Jdue y_n\in F$. Being $F$ a filter, from the last fact it follows that also $\Box\Jdue y_1\land...\land\Box\Jdue y_n\in F$, which implies $\Box\Jdue y_1\land...\land\Box\Jdue y_n=\Box(\Jdue y_1\land...\land\Jdue  y_n)=\Box z\in F$, obtained using distributivity of box over meet (K$_\land$). Since $\Box z\leq\Box\Jdue x$, $\Box\Jdue x\in F$ as well.
	
	Next we show that $\Jdue x\lor\Jzero x\in F$. Suppose not, therefore $\neg(\Jdue x\lor\Jzero x)=\neg\Jdue x\land\neg\Jzero x\in F$ (recall that these are all elements on the bottom fiber and there $F$ behaves like a standard ultrafilter), hence $\neg\Jdue x\land\neg\Jzero x\in F$, therefore $\neg\Jdue x,\neg\Jzero x\in F$. By definition of $b$-ultrafilter, the latter implies $x\notin F$ and $\neg x\notin F$, but since $F\cap A_x\neq\emptyset$ this contradicts Lemma \ref{lem.:jultrafilter_properties}, which requires that either $x\in F$ or $\neg x\in F$. Therefore $\Jdue x\lor\Jzero x\in F$. Now notice that $(\Box\Jdue x\land(\Jdue x\lor\Jzero x))\to\Jdue\Box x$ is a theorem of $\MBeK$ (hence of $\MBeKglobal$ as well), in fact take an arbitrary model $\mathfrak{M}=\pair{W,R,v}$ in which $v(w,\Box\Jdue x)=v(w,\Jdue x\lor\Jzero x)=1$ for $w\in W$. $v(w,\Box\Jdue x)=1$ implies that for all $w'\in R[w],v(w',\Jdue x)=v(w',x)=1$. On the other hand $v(w,\Jdue x\lor\Jzero x)=1$ requires either $v(w,\Jdue x)=1$ or $v(w,\Jzero x)=1$, which implies, respectively, $v(w,x)=1$ or $v(w,x)=0$. In both cases, $v(w,\Box x)=1$, hence $v(w,\Jdue\Box x)=1$.
	
	By algebraizability this theorem is equivalent to the validity of $(\Box\Jdue x\land(\Jdue x\lor\Jzero x))\to\Jdue\Box x\approx 1$ in $\MBCABoch$, and since the left-hand member of the equality is an element of the bottom fiber, this amounts to $\Box\Jdue x\land(\Jdue x\lor\Jzero x)\leq\Jdue\Box x$. As a consequence $\Jdue\Box x\in F$, hence $\Box x\in F$, contradicting the assumption. 
	
	This proves that the Boolean filter generated by $H'$ is non-trivial, hence it can be extended to an ultrafilter $H''$. Notice that $\neg\Jdue x\in H''$, hence $\Jdue x\notin H''$, therefore $x\notin H''$.
	
	Now extend $H''$ to a $b$-ultrafilter $U$, that is define 
	$$U=\bigcup_{y\in H''} \Jdue^{-1}[y].$$
	
	\noindent Clearly this is a $b$-ultrafilter extending $H$ and not containing $x$, furthermore by definition of $H$, $FQ_AU$, in fact if $\Box y\in F$, $\Jdue\Box y\in F$ and, by the disequation $\Jdue\Box y\leq\Jdue\Box\Jdue y$ valid in $\MBCABoch$ (which can be checked through algebraizability via the validity of $\Jdue\Box y\to\Jdue\Box\Jdue y$ on Bochvar-Kripke models), also $\Jdue\Box\Jdue y\in F$, hence $\Box\Jdue y\in F$, which implies $\Jdue y\in H\subseteq U$, then $y\in U$, therefore the two $b$-ultrafilters are related by $Q_A$. But then there exists $U\in Q_A[F]$ s.t. $x\notin U$, contradicting $Q_A[F]\subseteq\{G\in \textit{Uf}_B(\A)|x\in G\}$. This concludes the proof that $\Box x\in F$, which yields $\iota(\Box x)(F)=1$.
	
	For the remaining cases it is enough to observe that $\iota(\Box x)(F)=\ant$ iff $F\cap A_x=\emptyset$ iff $m_{Q_A}\iota(x)(F)=\ant$. This concludes the proof that $\iota$ is a homomorphism, hence it is the desired embedding.	
\end{proof}

In conclusion, we can employ the representation result of Theorem \ref{th.:MBCA_representation} to find an algebraic semantics for the local modal logic $\MBeK$.

\begin{lemma}\label{lem.:key_lemma_complex_frame}
	Let $\M=\pair{W,R,v}$ be a Bochvar-Kripke model and put $\A=\mathbf{Cm}_B(\pair{W,R})$. Let $h:\Fm\to\A$ be the homomorphism defined by $h(x)(w)=v(w,x)$, for all variables $x$ and $w\in W$. For all $\phi\in Fm$ and $w\in W,h(\phi)(w)=v(w,\phi)$.
\end{lemma}
\begin{proof}
	The proof is an easy induction over the complexity of formulae. The base case is provided by the definition. We then show just the case for box. $h(\Box\psi)(w)=(\Box h(\psi))(w)=1$ iff $h(\psi)(w)\neq\ant$ and $R[w]\subseteq h(\psi)^{-1}[\{1\}]=\{w'\in W|h(\psi)(w')=1\}$ iff, by induction hypothesis, $v(w,\psi)\neq\ant$ and for all $R[w]\subseteq\{w'\in W|v(w',\psi)=1\}$ iff $v(w,\Box\psi)=1$. Similarly $h(\Box\psi)(w)=(\Box h(\psi))(w)=\ant$ iff $h(\psi)(w)=\ant$ iff, by induction hypothesis, $v(w,\psi)=\ant$ iff $v(w,\Box\psi)=\ant$. The remaining cases are similar.
\end{proof}

\begin{theorem}\label{th.:local_logic_of_order}
	The local modal logic $\MBeK$ is the logic induced by the class of matrices  $\{\pair{ \A,U}$ $|$ $\A\in\MBCABoch,U\in\textit{Uf}_B(\A)\}$.
\end{theorem}
\begin{proof}
	For soundness, let $\A\in\MBCABoch$ and $U\in\textit{Uf}_B(\A)$. First notice that by algebraizability we have that every axiom $\phi$ of $\MBeK$ amounts to the validity of $\phi\approx 1$ in $\MBCABoch$; since $1\in U$, this proves the validity of all axioms. For the non-modal part of $\MBeK$ we can take the Hilbert-style system of \cite{FinnGrigolia} for $\Be$, whose only inference rule is modus ponens: given a valuation $h$, suppose $h(x),h(x\to y)\in U$, which holds iff $\Jdue h(x),\Jdue h(x\to y)=\Jdue (h(x)\to h(y))\in U$ by definition of $b$-ultrafilter, hence from the validity of $\Jdue\phi\land\Jdue(\phi\to\psi)\to\Jdue\psi\approx 1$ in $\BCA$ it follows $\Jdue h(x)\land\Jdue(h(x)\to h(y))\leq\Jdue h(y)$, so $h(y)\in U$. The sole remaining rule is the modal one, (local) necessitation: if $\phi$ is a theorem of $\MBeK$ again $\phi\approx 1$ is valid in $\MBCABoch$ by algebraizability, therefore $\Box\phi\approx\Box 1\approx 1$ holds as well, we conclude that $\Box\phi$ is a theorem.
	
	Adequacy will be proved by contraposition, by showing that every Bochvar-Kripke countermodel for a certain inference generates a modal Bochvar algebra that falsifies that same inference when equipped with an appropriate external ultrafilter. Let $\Gamma\nvdash_{\MBeK}\phi$. By completeness Theorem \ref{th: completezza Bochvar modale}, there is a Bochvar-Kripke countermodel $\M=\pair{W,R,v}$ where, for some $a\in W, v(a,\gamma)=1$ for every $\gamma\in\Gamma$, but $v(a,\phi)\neq 1$.
	
	Consider the corresponding full complex Bochvar algebra $\mathbf{Cm}_B(\pair{W,R})$, which by Lemma \ref{lem.:complex_alg_is_bochvar} belongs to $\MBCABoch$. It is easy to realize that the elements belonging to the bottom fiber of $\mathbf{Cm}_B(\pair{W,R})$ are exactly those belonging to the set $\{1,0\}^W$ (in fact an element $x$ of the bottom fiber on any $\BCA$ can be characterized as s.t. $x=\Jdue x$). Consider the map $f_a\in\{1,0\}^W$ defined, for all $w\in W$, as $f_a(w)=1$ iff $w=a$. $f$ is an atom of the bottom fiber, in fact it is a cover of 0, therefore the principal filter generated by $f_a$ in $\{1,0\}^W$, namely 
	$\uparrow f_a=\{g\in\{1,0\}^W$ $|$ $g(a)=1\}$, is an ultrafilter. 
	
	Expand $\uparrow f_a$ to a $b$-ultrafilter by defining $U_a=\bigcup_{g\in\uparrow f_a} \Jdue^{-1}[g]$. Next, consider the evaluation $h:\Fm\to\mathbf{Cm}_B(\pair{W,R})$ described in Lemma \ref{lem.:key_lemma_complex_frame}. By the mentioned lemma, we know that for all $\psi\in Fm$ and $w\in W,h(\psi)(w)=v(w,\psi)$, in particular for all $\gamma\in\Gamma,h(\gamma)(a)=1$, therefore $h(\gamma)\in U_a$. At the same time, $h(\phi)\notin U_a$, in fact $v(a,\phi)\neq 1$ iff $v(a,\Jdue\phi)\neq 1$ iff $h(\Jdue\phi)(a)\neq 1$ iff $h(\Jdue\phi)\notin\uparrow f_a$ iff $h(\Jdue\phi)=\Jdue h(\phi)\notin U_a$ iff $h(\phi)\notin U_a$. The matrix $\pair{\mathbf{Cm}_B(\pair{W,R}),U_a}$ is the desired countermodel.
\end{proof}

\subsection{Representation of modal Halld\'{e}n algebras}

Like we did in the previous case, we start from providing an appropriate notion of ultrafilter for $\MBCAHalld$.

\begin{definition}
	Let $\A$ be an algebra with a $\BCA$-reduct. $U\subseteq A$ is an \textit{paraconsistent external ultrafilter} ($h$-ultrafilter) if $\neg\Jzero[U]$ is a Boolean ultrafilter of the bottom fiber $\A_0$ and, for all $x\in A$, 
	$$x\in U\Leftrightarrow\neg\Jzero x\in U.$$ 
	The set of all $h$-ultrafilters over $\A$ will be denoted by $\textit{Uf}_H(\A)$.
\end{definition}

We use $h$-ultrafilters to build a Kripke frame out of an algebra.

\begin{definition}\label{def.:ultrafilter_frame_PWK}
	Given $\A\in\MBCAHalld$, its $h$-ultrafilter frame is the Kripke frame $\mathfrak{Uf}_H(\A)=\pair{\textit{Uf}_H(\A),S_A}$, with the accessibility relation defined as
	\begin{align*}
		FS_AG\Leftrightarrow\text{ for all }x\in A\text{ s.t. }F\cap A_x\neq A_x,\Box x\in F\text{ implies }x\in G.
	\end{align*}
%	\begin{align*}
%		FS_AG\Leftrightarrow\text{ for all }x\in A\text{ s.t. }F\cap A_x\neq A_x,\Box x\in F\text{ implies either }x\in G\text{ or }0_x\in G.
%	\end{align*}
\end{definition}

\begin{lemma}\label{lem.:hultrafilter_basics}
	For any $\A\in\BCA, U\in \textit{Uf}_H(\A)$ and $x,y\in A_i$ for some index $i\in I$ of its P\l onka sum representation. The following hold:
	\begin{itemize}
		\item if $x,y\in U$ then $x\land y\in U$;
		\item if $x\in U$ and $x\leq y$ then $y\in U$.
	\end{itemize}
\end{lemma}
\begin{proof}
	Let $x,y\in U$, then $\neg\Jzero x,\neg\Jzero y\in\neg\Jzero[U]$, which is an ultrafilter, hence $\neg\Jzero x\land\neg\Jzero y\in\neg\Jzero[U]$; observe that an immediate consequence of the P\l onka sum representation of $\BCA$ (Proposition \ref{prop_bca_isbl}) is that the fact that $x$ and $y$ belong to the same fiber can be expressed by $0_x=0_y$, hence by the quasiequation $0_x\approx 0_y\Rightarrow\neg\Jzero x\land\neg\Jzero y\approx\neg\Jzero(x\land y)$ valid in $\BCA$ (as can be checked on $\WKt$), we conclude $\neg\Jzero(x\land y)\in\neg\Jzero[U]$ and $x\land y\in U$.
	
	Now suppose that $x\leq y$. By Theorem \ref{th.: BCA_representation} $\Jdue$ is a Boolean isomorphism, in particular it is order-preserving, hence the composition $\neg\Jzero:=\neg\Jdue\neg$ is order-preserving as well, therefore $\neg\Jzero x\leq\neg\Jzero y$. If $x\in U$ then $\neg\Jzero x\in \neg\Jzero[U]$, then, since $\neg\Jzero[U]$ is an ultrafilter, $\neg\Jzero y\in \neg\Jzero[U]$ and $y\in U$.
\end{proof}

As we did in the previous representation theorem with $b$-ultrafilters, the current proof will rely on some properties of $h$-ultrafilters.

\begin{lemma}\label{lem.:hultrafilter_properties}
	Let $\A\in\BCA$ and $\mathbf{I}$ the semilattice of indexes of its P\l onka sum representation. For all $U\in \textit{Uf}_H(\A), i\in I$ and $x,y\in A$ the following hold:
	\begin{itemize}
		\item[(a)] $U\cap A_0\neq\emptyset$;
		\item[(b)] $1_i\in U$;
		\item[(c)] if $0_i\in U,U\cap A_i=A_i$;
		\item[(d)] if $U\cap A_i\neq A_i$, for all $x\in A_i$ exactly one of $x\in U$ or $\neg x\in U$ holds;
		\item[(e)] if $U\cap A_{i\lor j}\neq A_{i\lor j}$ then $U\cap A_{i}\neq A_{i}$;
		\item[(f)] if $U\cap A_{i}\neq A_{i}$ and $U\cap A_{j}\neq A_j$ then $U\cap A_{i\lor j}\neq A_{i\lor j}$;
	\end{itemize}
\end{lemma}
\begin{proof}
	(a) is an immediate consequence of the fact that $\neg\Jzero[U]$ is an ultrafilter over $\A_0$.
	
	For (b), $\neg\Jzero 1_i=\neg\Jzero(x\lor\neg x)=\neg\Jdue\neg(x\lor\neg x)=\neg\Jdue(x\land\neg x)=\neg\Jdue 0_x=1\in U$, again because $\neg\Jzero[U]$ is an ultrafilter.
	
	Point (c) is an immediate consequence of the upward closure of $U$ stated in Lemma \ref{lem.:hultrafilter_basics}.
	
	For (d), let $x\in U\cap A_i\neq A_i$. By point (b) we know that $1_i\in U$, so $\neg\Jzero 1_i\in U$, but $\neg\Jzero 1_i=\neg\Jzero(x\lor\neg x)=\neg\Jzero x\lor\neg\Jzero\neg x$ and $U$ is an ultrafilter over $\A_0$, therefore at least one of $\neg\Jzero x\in U$ or $\neg\Jzero\neg x\in U$ holds, hence $x\in U$ or $\neg x\in U$. That exactly one holds follows from the fact that otherwise $x\land\neg x=0_x\in U$, contradicting the hypothesis by point (c).

	For (e), suppose $U\cap A_{i\lor j}\neq A_{i\lor j}$, then by upward closure of $U$ it must hold $0_{i\lor j}\notin U$, so $\neg\Jzero 0_{i\lor j}\notin U$. Since $\neg\Jzero 0_{i\lor j}\geq\neg\Jzero 0_{i}$ (as it can be checked on $\WKt$), again by upward closure $\neg\Jzero 0_{i}\notin U$, hence $0_{i}\notin U$.
	
	For (f), the assumptions $U\cap A_{i}\neq A_{i}$ and $U\cap A_{j}\neq A_j$ imply $0_i,0_j\notin U$, using point (c) by contraposition. Hence $\neg\Jzero 0_i,\neg\Jzero 0_j\notin U$, then since $U$ is an ultrafilter it follows $\neg\Jzero 0_i\lor\neg\Jzero 0_j=\neg\Jzero(0_i\lor 0_j)\notin U$, therefore $0_i\lor 0_j\notin U$, which proves the claim since $0_i\lor 0_j\in A_{i\lor j}$ by the computation of operations in a P\l onka sum.
\end{proof}

\noindent The structure of $h$-ultrafilters is quite different from that of $b$-ultrafilters, in fact the former contain all local tops and, in particular, the fix-point for negation. Moreover Lemma \ref{lem.:hultrafilter_basics} together with points (c) and (d) of Lemma \ref{lem.:hultrafilter_properties} yields the following corollary.

\begin{corollary}\label{cor.:hultrafilters_boolean}
	Let $\A\in\BCA$ and $U\in \textit{Uf}_H(\A)$. If $U\cap A_i\neq A_i$ then $U\cap A_i$ is an ultrafilter over $\A_i$.
\end{corollary}

\n Therefore every fiber of a $\MBCAHalld$ is either completely or partially covered by any $h$-ultrafilter, and in the former case that $h$-ultrafilter is a Boolean ultrafilter over the fiber. 

\begin{figure}[h]
	\centering
	\begin{tikzpicture}[scale=1]
		\node (nn) at (0,6.5) {$(\ant,\ant)$};
		\node (1n) at (-3,5) {$(1,\ant)$};
		\node (0n) at (-3,3) {$(0,\ant)$};
		\node (n1) at (3,5) {$(\ant,1)$};
		\node (n0) at (3,3) {$(\ant,0)$};
		\node (11) at (0,3) {$(1,1)$};
		\node (01) at (-1.5,1.5) {$(0,1)$};
		\node (10) at (1.5,1.5) {$(1,0)$};
		\node (00) at (0,0) {$(0,0)$};
		
		\draw (00) -- (01) -- (11) -- (10) -- (00);
		\draw (0n) -- (1n);
		\draw (n0) -- (n1);
		
		\draw (-0.75,2.25) circle [x radius=0.8cm, y radius=1.7cm, rotate=-45];
		\draw[dashed] (0.75,2.25) circle [x radius=0.8cm, y radius=1.7cm, rotate=45];
		\draw (-3,4) circle [x radius=0.9cm, y radius=1.75cm];
		\draw[dashed] (3,4) circle [x radius=0.9cm, y radius=1.75cm];
		\draw (3,5) circle [radius=0.6cm];
		\draw[dashed] (-3,5) circle [radius=0.6cm];
		\draw (0,6.5) circle [radius=0.85cm];
		\draw[dashed] (0,6.5) circle [radius=0.75cm];
		
		\draw [->,dotted] (1n) to [out=0,in=120,looseness=1] (11);
		\draw [->,dotted] (0n) to [out=-90,in=180,looseness=1] (01);
		\draw [->,dotted] (n1) to [out=180,in=60,looseness=1] (11);
		\draw [->,dotted] (n0) to [out=-90,in=0,looseness=1] (10);
		\draw [->,dotted] (nn) to [out=-90,in=90] (11);
	\end{tikzpicture}
	\caption{$\WKt^2$ with its two $h$-ultrafilters (dotted lines represent $\neg\Jzero$).}\label{fig.: WKt_h_ultrafilter}
\end{figure}
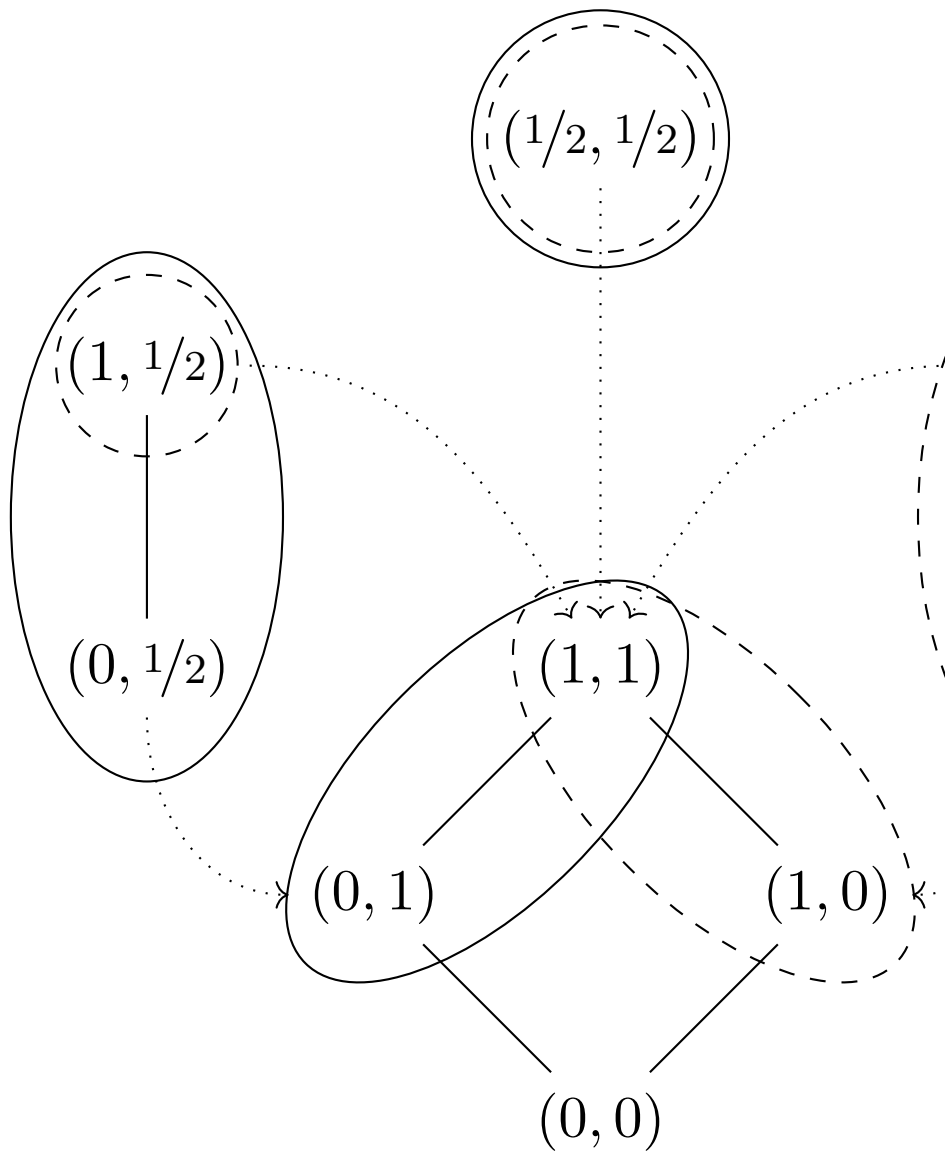

From the previous facts it is easy to conclude that for $\A\in\MBCAHalld$ there is a bijective correspondence between its external ultrafilters and the Boolean ultrafilters of its bottom fiber $\A_0$. As an example, Figure \ref{fig.: WKt_h_ultrafilter} shows the two $h$-ultrafilters of $\WKt^2$ (represented according to its P\l onka sum decomposition).

\begin{definition}\label{def.:complex_hallden_alg}
	Given a Kripke frame $\mathfrak{F}=\pair{W,R}$, its \textit{full Halld\'en complex algebra} is $\mathbf{Cm}_H(\mathfrak{F})=\pair{\{1,\ant,0\}^W,\lor,\neg,\Jdue,n_R,1,0}$, which differs from a complex Bochvar algebra of Definition \ref{def.:bochvar_complex_alg} only for
	$$
	(n_R f)(w)=\begin{cases} 
			1 & \text{if }f(w)\neq\ant\text{ and }R[w]\subseteq f^{-1}[\{1,\ant\}];\\
			\ant & \text{if }f(w)=\ant;\\
			0 & \text{otherwise.}	
		\end{cases}
	$$
	
	A \textit{Halld\'en complex algebra} of $\mathfrak{F}$ is any subalgebra of $\mathbf{Cm}_H(\mathfrak{F})$.	
\end{definition}

\begin{lemma}\label{lem.:complex_alg_is_hallden}
	For any Kripke frame $\mathfrak{F},\mathbf{Cm}_H(\mathfrak{F})\in\MBCAHalld$.
\end{lemma}
\begin{proof}
	Let $\A=\mathbf{Cm}_H(\mathfrak{F})$ for some Kripke frame $\mathfrak{F}=\pair{W,R}$. That the $\Box$-free reduct of $\A$ is a Bochvar algebra is immediate by its definition, since $\A$ is isomorphic to the power $\WKt^{|W|}$. Next we prove that $\A$ satisfies the modal axioms of $\MBCAHalld$. Since all the considered equations involve only external formulae, the cases for $\ant$ will be discarded.
	
	\vspace{0.2cm}
	
	(M1) The proof proceeds essentially as the same case in Lemma \ref{lem.:complex_alg_is_bochvar}.
	
	\vspace{0.2cm}
	
	(M2) $(\neg\Juno f)(w)=1$ iff $f(w)\neq\ant$ iff $(\neg\Juno\Box f)(w)=1$.
	
	\vspace{0.2cm}
	
	(MP3) By contradiction, suppose $\neg\Jzero(\Box x\leftrightarrow\Box\neg\Jzero x)(w)=0$ for some $w\in W$. This holds iff $(\Box x\leftrightarrow\Box\neg\Jzero x)(w)=0$. Now in the case that $x(w)=\ant$ then $(\Box x\leftrightarrow\Box\neg\Jzero x)(w)=\ant$, contradicting the assumption. On the case that $x(w)\neq\ant$, $(\Box x)(w)=1$ iff $R[w]\subseteq x^{-1}[\{1,\ant\}]=\Jzero x^{-1}[\{0\}]=\neg\Jzero x^{-1}[\{1\}]$ iff $(\Box\neg\Jzero x)(w)=1$ (since always $(\neg\Jzero x)(w)\neq\ant$), therefore $(\Box x\leftrightarrow\Box\neg\Jzero x)(w)=1$, so in both cases we conclude a contradiction.
	
	\vspace{0.2cm}
	
	($\tau$-GN$'$) Suppose for all $w\in W,(\neg\Jzero x)(w)=1$, that is $x(w)\neq 0$. So $R[w]\subseteq W=x^{-1}[\{1,\ant\}]$, hence in the case that $x(w)\neq\ant$ then $(\Box x)(w)=1=(\neg\Jzero\Box x)(w)$, and equally in the other case $x(w)=\ant$ implies $(\neg\Jzero\Box x)(w)=1$.
	
\end{proof}

\begin{theorem}[Representation]\label{th.:MBCAHalld_representation}
	Let $\A\in\MBCAHalld$. There exists an embedding $\eta:\A\to\mathbf{Cm}_H(\mathfrak{Uf}_H(\A))$ given by, for all $x\in A$ and $F\in \textit{Uf}_H(\A)$:
	$$
	\eta(x)(F)=\begin{cases} 
			1 & \text{if }F\cap A_x\neq A_x\text{ and }x\in F;\\
			\ant & \text{if }F\cap A_x= A_x;\\
			0 & \text{otherwise.}
		\end{cases}
	$$
\end{theorem}
\begin{proof}
	To prove that the map $\eta$ is injective we consider for $x\neq y\in A$ two cases. If $\neg\Jzero x\neq\neg\Jzero y$ then there exists an ultrafilter $U$ over $\A_0$ which separates the two points, that is, w.l.o.g., $\neg\Jzero x\in U$ and $\neg\Jzero y\notin U$. We can expand $U$ to the $h$-ultrafilter $V=\bigcup_{z\in U} \neg\Jzero^{-1}[z]$, for which $x\in V$ and $y\notin V$. Hence $\eta(y)(V)=0\neq\eta(x)(V)$, so $\eta(x)\neq\eta(y)$. If $\neg\Jzero x=\neg\Jzero y$, it must be the case that $A_x\neq A_y$, because if $x$ and $y$ belonged to the same fiber their images through $\neg\Jzero:=\neg\Jdue\neg$ cannot coincide, since $\Jdue$ restricted to that fiber is an isomorphism. Since $A_x\neq A_y$, also $\Jzero 0_x=\Jdue 1_x\neq\Jdue 1_y=\Jzero 0_y$ (Theorem \ref{th.: BCA_representation}). Now consider an ultrafilter $U$ s.t. $\Jzero 0_x\in U$ and $\Jzero 0_y\notin U$ and extend it to the $h$-ultrafilter $V$ defined in the previous case. $\Jzero 0_x\in V$ and $\Jzero 0_y\notin V$ hold iff $\neg\Jzero 0_x\notin V$ and $\neg\Jzero 0_y\in V$ iff $0_x\notin V$ and $0_y\in V$. Since $V$ is upward closed (Lemma \ref{lem.:hultrafilter_basics}), $0_y\in V$ implies $V\cap A_y=A_y$. We conclude $\eta(x)(V)\neq\ant=\eta(y)(V)$.
	
	Now we show that $\eta$ is a homorphism.
	
	\vspace{0.2cm}
	
	($\lor$) Suppose $\eta(x\lor y)(F)=1$, which holds iff $F\cap A_{x\lor y}\neq A_{x\lor y}$ and $x\lor y\in F$. The latter implies $\neg\Jzero(x\lor y)=\neg\Jzero x\lor\neg\Jzero y\in F$ (easy to check on $\WKt$). Since $F$ is an ultrafilter over $\A_0$, $\neg\Jzero x\lor\neg\Jzero y\in F$ implies $\neg\Jzero x\in F$ or $\neg\Jzero y\in F$, which holds iff $x\in F$ or $y\in F$. By Lemma \ref{lem.:hultrafilter_properties}, point (e), $F\cap A_{x\lor y}\neq A_{x\lor y}$ implies $F\cap A_{x}\neq A_{x}$ and $F\cap A_{y}\neq A_{y}$, therefore either $\eta(x)(F)=1$ and $\eta(y)(F)\neq\ant$ or $\eta(y)(F)=1$ and $\eta(x)(F)\neq\ant$. Hence $(\eta(x)\lor\eta(y))(F)=1$.
	
	Next, suppose $(\eta(x)\lor\eta(y))(F)=1$, that is $\eta(x)(F)=1$ and $\eta(y)(F)\neq\ant$, or $\eta(x)(F)\neq\ant$ and $\eta(y)(F)=1$. W.l.o.g., let us assume the first disjunct, which is equivalent to $x\in F$ and $F\cap A_y\neq A_y$, then by Lemma \ref{lem.:hultrafilter_properties}, point (f) we obtain $F\cap A_{x\lor y}\neq A_{x\lor y}$. Now $x\in F$ iff $\neg\Jzero x\in F$, moreover $\neg\Jzero x\leq\neg\Jzero(x\lor y)$, therefore $x\lor y\in F$. We conclude that $\eta(x\lor y)(F)=1$.
	
	Finally assume $\eta(x\lor y)(F)=\ant$, which holds iff $F\cap A_{x\lor y}=A_{x\lor y}$. Using Lemma \ref{lem.:hultrafilter_properties}, points (e)-(f), the latter holds iff either $F\cap A_x=A_x$ or $F\cap A_y=A_y$ iff $\eta(x)(F)=\ant$ or $\eta(y)(F)=\ant$ iff $(\eta(x)\lor\eta(y))(F)=\ant$.

	\vspace{0.2cm}
	
	($\neg$) Suppose $\eta(\neg x)(F)=1$, which holds iff $F\cap A_{\neg x}=F\cap A_{x}\neq A_x$ and $\neg x\in F$ iff $F\cap A_{x}\neq A_x$ and $x\notin F$ (by Lemma \ref{lem.:hultrafilter_properties}) iff $\eta(x)(F)=0$ iff $(\neg\eta(x))(F)=1$. Next assume $\eta(\neg x)(F)=\ant$, which holds iff $F\cap A_{x}=A_{x}$ iff $\eta(x)(F)=\ant$ iff $(\neg\eta(x))(F)=\ant$.
	
	\vspace{0.2cm}
	
	($\Jdue$) Suppose $\eta(\Jdue x)(F)=1$, which holds iff $F\cap A_{\Jdue x}\neq A_{\Jdue x}$ (which always holds since $A_{\Jdue x}=A_0$) and $\Jdue x\in F$. Now since $F$ is an $h$-ultrafilter, $\Jzero\neg x\in F$ iff $\neg\Jzero\neg x\notin F$ iff $\neg x\notin F$ iff $x\in F$ and $F\cap A_x\neq A_x$. The last two facts are equivalent to $(\Jdue\eta(x))(F)=1$. Since $(\Jdue\eta(x))(F)$ and $\eta(\Jdue x)(F)$ can never return $\ant$, this is enough to prove the case.
	
	\vspace{0.2cm}
	
	($\Box$) We want to prove $\eta(\Box x)=n_{S_A}\eta(x)$. Let us first make explicit the clauses of $n_{S_A}\eta(x)$. For $F\in \textit{Uf}_H(\A)$,
	$$
	(n_{S_A}\eta(x))(F)=\begin{cases} 
			1 & \text{if }\eta(x)(F)\neq\ant\text{ and }S_A[F]\subseteq \eta(x)^{-1}[\{1,\ant\}];\\
			\ant & \text{if }\eta(x)(F)=\ant;\\
			0 & \text{otherwise}	
		\end{cases}
	$$
	which can be equivalently rephrased as
	$$
	(n_{S_A}\eta(x))(F)=\begin{cases} 
			1 & \text{if }F\cap A_x\neq A_x\text{ and }S_A[F]\cap\{G\in \textit{Uf}_H(\A)\text{ | }G\cap A_x\neq A_x\\
			& \text{ and }\neg x\in G\}=\emptyset;\\
			\ant & \text{if }F\cap A_x=A_x;\\
			0 & \text{otherwise.}	
		\end{cases}
	$$
	
	Suppose that $\eta(\Box x)(F)=1$, therefore $F\cap A_{\Box x}\neq A_x$ and $\Box x\in F$. By Lemma \ref{lem.: MBCA_box_in_fibra}, we know that every fiber of the P\l onka representation of $\A$ is closed under $\Box$, therefore $A_{\Box x}=A_x$, hence $F\cap A_x\neq A_x$. Since $\Box x\in F$, Definition \ref{def.:ultrafilter_frame_PWK} implies that for all $G\in S_A[F]$, either $x\in G$ or $0_x\in G$. If $0_x\in G,G\cap A_x=A_x$ (Lemma \ref{lem.:hultrafilter_properties}), then $\eta(x)(G)=\ant$, while if $x\in G$ then $\eta(x)(G)\neq 0$, therefore $S_A[F]\subseteq \eta(x)^{-1}[\{1,\ant\}]$. This amounts to $(n_{S_A}\eta(x))(F)=1$.
	
	Next suppose that $(n_{S_A}\eta(x))(F)=1$, so $S_A[F]\cap\{G\in \textit{Uf}_H(\A)\text{ | }G\cap A_x\neq A_x$ and $\neg x\in G\}=\emptyset$ and $F\cap A_x\neq A_x$. Since $A_x=A_{\Box x}$, all we are left to show in order to prove $\eta(\Box x)(F)=1$ is $\Box x\in F$. Assume the contrary, that is $\Box x\notin F$, and consider the sets
	$$
	H=\{\neg\Jzero y\text{ | }\Box\neg\Jzero y\in F\},
	$$$$
	H'=H\cup\{\Jzero x\}.
	$$
	Observe that $H'$ belongs entirely to the bottom fiber $\A_0$ of the P\l onka representation. To demonstrate that the Boolean filter generated by $H'$ is non-trivial, we need to prove that no finite meet of elements of $H'$ is equal to 0. To do that take $z$ finite meet of elements of $H$ and suppose towards contradiction that $z\land\Jzero x=0$. Using basic properties of Boolean algebras this implies $z\leq\neg\Jzero x$, then $\Box z\leq\Box\neg\Jzero x$ by monotonicity of $\Box$ (recall that $\A_0$ is a modal algebra, Lemma \ref{lem.: bottom_bca_bao}). By definition $z$ is of the form $\neg\Jzero y_1\land...\land\neg\Jzero y_n$ for $\neg\Jzero y_1,...\neg\Jzero y_n\in H$, therefore $\Box\neg\Jzero y_1,...\Box\neg\Jzero y_n\in F$. Being $F$ a filter, from the last fact it follows that also $\Box\neg\Jzero y_1\land...\land\Box\neg\Jzero y_n\in F$, which implies $\Box\neg\Jzero y_1\land...\land\Box\neg\Jzero y_n=\Box(\neg\Jzero y_1\land...\land\neg\Jzero  y_n)=\Box z\in F$, obtained using distributivity of box over meet (K$_\land$) (again, because $\A_0$ is a modal algebra). Since $\Box z\leq\Box\neg\Jzero x$, $\Box\neg\Jzero x\in F$ as well.
	
	Now notice that $\Box\neg\Jzero x\to\neg\Jzero\Box x$ is a theorem of $\MPWK$ (hence of $\MPWKglobal$ as well), in fact take an arbitrary model $\mathfrak{M}=\pair{W,R,v}$ in which $v(w,\Box\neg\Jzero x)=1$ but $v(w,\neg\Jzero\Box x)=0$. The latter implies $v(w,\Jzero\Box x)=1$, hence $v(w,\Box x)=0$, that is $v(w,x)\neq\ant$ and there exists $w'\in R[w]$ such that $v(w',x)=0$. By $v(w,\Box\neg\Jzero x)=1$ we have $v(w',\neg\Jzero x)\neq 0$, which gives $v(w',\neg\Jzero x)=1$, that holds if either $v(w',x)=1$ or $v(w',x)=\ant$, both contradicting $v(w',x)=0$. Therefore no such countermodel exists. By algebraizability this theorem is equivalent to the validity of $\Box\neg\Jzero x\to\neg\Jzero\Box x\approx 1$ in $\MBCAHalld$, and since the left-hand member of the equality is an element of the bottom fiber, this amounts to $\Box\neg\Jzero x\leq\neg\Jzero\Box x$. As a consequence $\neg\Jzero\Box x\in F$, hence $\Box x\in F$, contradicting the assumption. 
	
	This proves that the Boolean filter generated by $H'$ is non-trivial, hence it can be extended to an ultrafilter $H''$. Notice that $\Jzero x\in H''$, hence $\neg\Jzero x\notin H''$ (Lemma \ref{lem.:hultrafilter_properties}), therefore $x\notin H''$.
	
	Now extend $H''$ to an $h$-ultrafilter $U$, that is define 
	$$U=\bigcup_{y\in H''} \neg\Jzero^{-1}[y].$$
	
	\noindent Clearly this is an $h$-ultrafilter extending $H$ and not containing $x$. We show that $FS_AU$. Let $\Box y\in F$ s.t. $F\cap A_y\neq A_y$, Lemma \ref{lem.:hultrafilter_properties} implies $\neg\Box y\notin F$, which holds iff $\neg\Jzero\neg\Box y\notin F$, therefore $\Jzero\neg\Box y\in F$ (because $F$ is a Boolean ultrafilter on the bottom fiber). By the disequation $\Jzero\neg\Box y\leq\Box\neg\Jzero y$, whose validity in $\MBCAHalld$ can be checked through algebraizability via the validity of $\Jzero\neg\Box y\to\Box\neg\Jzero y$ on Halld\'{e}n-Kripke models, we have $\Box\neg\Jzero y\in F$, which by definition of $H$ implies $\neg\Jzero y\in H\subseteq U$, therefore $y\in U$.
	
	We have proved that the two $h$-ultrafilters are related by $S_A$, but this means there exists $U\in S_A[F]$ s.t. $x\notin U$, that is $G\cap A_x\neq A_x$ and $\neg x\in G$, contradicting $S_A[F]\cap\{G\in \textit{Uf}_H(\A)\text{ | }G\cap A_x\neq A_x$ and $\neg x\in G\}=\emptyset$. This concludes the proof that $\Box x\in F$, which yields $\eta(\Box x)(F)=1$.
	
	For the remaining cases it is enough to observe that $\eta(\Box x)(F)=\ant$ iff $F\cap A_x=A_x$ iff $n_{S_A}\eta(x)(F)=\ant$. This concludes the proof that $\eta$ is a homomorphism, hence it is the desired embedding.
\end{proof}

With this representation we can easily adapt the proof Theorem \ref{th.:local_logic_of_order} to find an algebraic semantics for the local modal logic $\MPWK$.

\begin{theorem}
	The local modal logic $\MPWK$ is the logic induced by the class of matrices  $\{\pair{ \A,U}$ $|$ $\A\in\MBCAHalld,U\in\textit{Uf}_H(\A)\}$.
\end{theorem}
\begin{proof}
	We adapt Theorem \ref{th.:local_logic_of_order}. For soundness, let $\A\in\MBCAHalld$ and $U\in\textit{Uf}_H(\A)$. By algebraizability, for every axiom $\phi$ of $\MPWK$ we have that $\neg\Jzero\phi\approx 1$ in $\MBCAHalld$; since $1\in U$, it follows $\neg\Jzero\phi\in U$ iff $\phi\in U$ by definition of $h$-ultrafilter. For the non-modal part of $\MPWK$ we can take the Hilbert-style system of \cite{Segerberg65} for $\PWKe$, whose only inference rule is the following restricted modus ponens: $\phi,\phi\to\psi\vdash\psi$, given that every non-convered variable of $\phi$ is non-convered in $\psi$ as well (a variable is non-convered when it doesn't fall under the scope of any $\Jdue$ operator). Now take a valuation $h$ and suppose $h(x),h(x\to y)\in U$, which holds iff $\neg\Jzero h(x),\neg\Jzero h(x\to y)=\neg\Jzero (h(x)\to h(y))\in U$. Given the constraint, $\neg\Jzero\phi\land\neg\Jzero(\phi\to\psi)\to\neg\Jzero\psi\approx 1$ is valid in $\WKt$, hence in $\BCA$, since it would fail only for $h(\psi)=0$ and $h(\phi)=\ant$, which is forbidden by the requirement on variables, in fact $h(\phi)=\ant$ iff there is a non-convered variable $z$ in $\phi$ s.t. $h(z)=\ant$, but then $z$ should belong to $\psi$ and so $h(\psi)=\ant\neq 0$. It follows $\neg\Jzero h(x)\land\neg\Jzero(h(x)\to h(y))\leq\neg\Jzero h(y)$, so $h(y)\in U$. We are left to check the rule of necessitation: if $\phi$ is a theorem of then $\neg\Jzero\phi\approx 1$ is valid in $\MBCAHalld$ by algebraizability, then by ($\tau$-GN$'$) $\neg\Jzero\Box\phi\approx 1$. Since $1\in U,\neg\Jzero\Box\phi\in U$ as well, which implies $\Box\phi\in U$. 
	
	Adequacy will be proved by contraposition, again by showing that every Halld\'en-Kripke countermodel for a certain inference generates a modal Halld\'en algebra that falsifies that same inference when equipped with an appropriate external ultrafilter. Let $\Gamma\nvdash_{\MPWK}\phi$. By completeness Theorem \ref{th: completezza MPWK}, there is a Halld\'{e}n-Kripke countermodel $\M=\pair{W,R,v}$ where, for some $a\in W, v(a,\gamma)\neq 0$ for every $\gamma\in\Gamma$, but $v(a,\phi)=0$.
	
	Consider the corresponding full complex Halld\'{e}n algebra $\mathbf{Cm}_H(\pair{W,R})$, which by Lemma \ref{lem.:complex_alg_is_hallden} belongs to $\MBCAHalld$. By the same reasoning of  Theorem \ref{th.:local_logic_of_order}, the bottom fiber of $\mathbf{Cm}_H(\pair{W,R})$ has $\{1,0\}^W$ as support. Consider the same map $f_a\in\{1,0\}^W$ defined in the previous theorem, that is $f_a(w)=1$ iff $w=a$, for all $w\in W$. As before, $f$ is an atom of the bottom fiber, in fact it is a cover of 0, therefore the principal filter generated by $f_a$ in $\{1,0\}^W$, namely 
	$\uparrow f_a=\{g\in\{1,0\}^W$ $|$ $g(a)=1\}$, is an ultrafilter. 
	
	Expand $\uparrow f_a$ to an $h$-ultrafilter by defining $U_a=\bigcup_{g\in\uparrow f_a} \neg\Jzero^{-1}[g]$. Next, it is immediate to adapt the evaluation $h:\Fm\to\mathbf{Cm}_H(\pair{W,R})$ described in Lemma \ref{lem.:key_lemma_complex_frame} to the case of Halld\'{e}n algebras and prove that the lemma still holds. Hence we know that for all $\psi\in Fm$ and $w\in W,h(\psi)(w)=v(w,\psi)$, in particular for all $\gamma\in\Gamma,h(\gamma)(a)\neq 0$ iff $h(\neg\Jzero\gamma)(a)=1$ iff $\neg\Jzero\gamma\in\uparrow f_a\subseteq U_a$ iff $\gamma\in U_a$. On the other hand, $v(a,\phi)=0$ iff $v(a,\neg\Jzero\phi)=h(\neg\Jzero\phi)(a)=0$ iff $h(\neg\Jzero\phi)\notin\uparrow f_a$ iff $h(\neg\Jzero\phi)=\neg\Jzero h(\phi)\notin U_a$ iff $h(\phi)\notin U_a$. The matrix $\pair{\mathbf{Cm}_H(\pair{W,R}),U_a}$ is the desired countermodel.
\end{proof}

%\begin{theorem}
%	The local modal logic $\MPWK$ is the logic of order $\vdash_{\MBCAHalld}^{\leq_\land}$.
%\end{theorem}
%\begin{proof}
%	It is straightforward adaptation of Theorem \ref{th.:local_logic_of_order}. For adequacy, completeness Theorem \ref{th: completezza MPWK} gives a Halld\'{e}n-Kripke countermodel $\M=\pair{W,R,v}$ for $\Gamma\nvdash_{\MPWK}\phi$, where for some $a\in W, v(a,\gamma)\neq 0$ for every $\gamma\in\Gamma$, but $v(a,\phi)=0$.
%	
%	We then consider the corresponding complex Halld\'{e}n algebra $\mathbf{Cm}_H(\pair{W,R})$ and define an evaluation $h$ as described in Lemma \ref{lem.:key_lemma_complex_frame}. It is immediate to prove that the lemma holds in the case of modal Halld\'{e}n algebras as well, therefore for all $\psi\in Fm$ and $w\in W,h(\psi)(w)=v(w,\psi)$. We consider the map $f\in\{1,\ant,0\}^W$ defined, for all $w\in W$, as $f(w)=\min\limits_{\gamma\in\Gamma} \{h(\gamma)(w)\}$, where min is the minimum function in the order of meet. Clearly for all $\gamma\in\Gamma,f\leq_{\land}h(\gamma)$ and in particular $f(a)\neq 0$, therefore $f(a)\nleq_{\land}h(\phi)(a)=0$. This proves $\Gamma\nvdash^{\leq_{\land}}_{\MBCAHalld}\phi$.
%\end{proof}

\subsection{Term-equivalence} We conclude the paper by showing how the two quasivarieties for modal weak Kleene logics are term-equivalent. We start by noticing how the modalities of the languages of $\MBeK$ and $\MPWK$ are interdefinable.

\begin{lemma}\label{lem.: kripke_models_translation}
	Let $\pair{W,R}$ be a Kripke frame, $v_B$ and $v_H$ valuations s.t. $\M_B=\pair{W,R,v_B}$ is a Bochvar-Kripke model, $\M_H=\pair{W,R,v_H}$ is a Halld\'en-Kripke model, moreover $v_B$ and $v_H$ coincide w.r.t. variables. Define the translation maps:
	\begin{itemize}
		\item $\mu:Fm\to Fm$ as $\mu(x)=x$ for all variables and $\mu$ is a homomorphism except for $\mu(\Box\phi)=\Box\Jdue\mu(\phi)\land(\Box\Jdue\mu(\phi)\lor\mu(\phi))$;
		\item $\nu:Fm\to Fm$ as $\nu(x)=x$ for all variables and $\nu$ is a homomorphism except for $\nu(\Box\phi)=\Box\neg\Jzero\nu(\phi)\land(\Box\neg\Jzero\nu(\phi)\lor\nu(\phi))$.
	\end{itemize}
	It holds for all $w\in W,v_B(w,\phi)=v_H(w,\mu(\phi))$ and $v_H(w,\phi)=v_B(w,\  \nu(\phi))$. Furthermore the two translations are mutually inverse, in the sense that $v_B(w,\phi)=v_B(w,\nu\mu(\phi))$ and $v_H(w,\phi)=v_H(w,\mu\nu(\phi))$.
\end{lemma}
\begin{proof}
	It is an easy induction over the complexity of formulae. We illustrate the case for $\Box$ when considering $\mu$. Suppose $v_H(w,\mu(\Box\phi))=v_H(w,\Box\Jdue\mu(\phi)\land(\Box\Jdue\mu(\phi)\lor\mu(\phi)))=1$, which holds iff $v_H(w,\Box\Jdue\mu(\phi))=1$ and $v_H(w,\Box\Jdue\mu(\phi)\lor\mu(\phi))=1$. The latter is equivalent to either $v_H(w,\Box\Jdue\mu(\phi))=1$ and $v_H(w,\mu(\phi))\neq\ant$ or $v_H(w,\mu(\phi))=1$. Since in both cases $v_H(w,\Box\Jdue\mu(\phi))=1$, the latter case is a subcase of the former, so we conclude that $v_H(w,\mu(\phi))\neq\ant$ and for all $w'\in R[w],v_H(w',\Jdue\mu(\phi))=1$, that is $v_H(w',\mu(\phi))=1$, iff, by induction hypothesis, $v_B(w,\phi)\neq\ant$ and for all $w'\in R[w],v_B(w',\phi)=1$, that is $v_B(w,\Box\phi)$. The case for $\ant$ is straightforward.
	
	The proof for $\nu$ runs similarly. That the translations are mutually inverse can be easily checked by showing how the semantic clauses coincide between $\phi$ and $\nu\mu(\phi)$ (resp., $\mu\nu(\phi)$) in Bochvar-Kripke models (in Halld\'en-Kripke model).
\end{proof}

\n Therefore any Bochvar-Kripke model is equivalent with an Halld\'en-Kripke model with the same frame and a valuation coinciding on variables, and viceversa. The term-equivalence of $\MBCABoch$ and $\MBCAHalld$ is stated by the next theorem.

\begin{theorem}
	Consider the translations of Lemma \ref{lem.: kripke_models_translation}. The following hold: 
	$$
	\Theta\vDash_{\MBCABoch}\alpha\approx\beta\Leftrightarrow\mu[\Theta]\vDash_{\MBCAHalld}\mu(\alpha\approx\beta),
	$$$$
	\Theta\vDash_{\MBCAHalld}\alpha\approx\beta\Leftrightarrow\nu[\Theta]\vDash_{\MBCABoch}\nu(\alpha\approx\beta).
	$$
	Furthermore the two translations are mutually inverse, in the sense that 
	$$
	\MBCABoch\vDash\alpha\approx\nu\mu(\alpha),
	$$$$
	\MBCAHalld\vDash\alpha\approx\mu\nu(\alpha).
	$$
\end{theorem}
\begin{proof}
	Let $\Theta\vDash_{\MBCABoch}\alpha\approx\beta$, which by (ALG2) of algebraizability holds iff $\{\epsilon\equiv\delta$ $|$ $\epsilon\approx\delta\in\Theta\}\vdash_{\MBeKglobal}\alpha\equiv\beta$ iff $\{\epsilon\equiv\delta$ $|$ $\epsilon\approx\delta\in\Theta\}\vDash^g_{\MBeK}\alpha\equiv\beta$, which is equivalent by Lemma \ref{lem.: kripke_models_translation} to $\{\mu(\epsilon)\equiv\mu(\delta)$ $|$ $\epsilon\approx\delta\in\Theta\}\vDash^g_{\MPWK}\mu(\alpha)\equiv\mu(\beta)$ iff $\{\mu(\epsilon)\equiv\mu(\delta)$ $|$ $\epsilon\approx\delta\in\Theta\}\vdash_{\MPWKglobal}\mu(\alpha)\equiv\mu(\beta)$ iff $\mu[\Theta]\vDash_{\MBCAHalld}\mu(\alpha)\approx\mu(\beta)$. That the translations are mutually inverse is a consequence of $\vdash_{\MBeKglobal}\alpha\equiv\nu\mu(\alpha)$ (an immediate corollary of Lemma \ref{lem.: kripke_models_translation}), hence by (ALG2) $\MBCABoch\vDash\alpha\approx\nu\mu(\alpha)$. It is simple to adapt the proof for the remaining claims.
\end{proof}

%\section{Conclusion}\label{sct:Conclusion}

\section*{Declaration on the use of AI tools}
\noindent The ideas, results and proofs presented in this paper are entirely due to the authors. Generative AI tools, namely Claude (Anthropic), were used to assist in checking the formal correctness of the arguments. All outputs of these tools were critically reviewed and verified by the authors, who take full responsibility for the content of this paper.

\section*{Acknowledgements}

\noindent All authors acknowledge the support of the Italian Ministry of University and Research through the PRIN 2022 project DeKLA (``Developing Kleene Logics and their Applications'', project code: 2022SM4XC8).

S. Bonzio acknowledges the support of Fondazione di Sardegna through the project LoMEA, the European Commission through the project SemPER (Grant agreement ID: 101299559) and of INdAM-GNSAGA (Gruppo Nazionale per le Strutture Algebriche, Geometriche e loro Applicazioni).

N. Zamperlin acknowledges the support of Fondazione di Sardegna through the project ``An algebraic approach to hyperintensionality'' (project code: F23C25000360007).

\bibliographystyle{abbrv}
\bibliography{references}

@book{Font16_handbook,
	author = {Josep Maria Font},
	publisher = {College Publications},
	title = {Abstract Algebraic Logic: An Introductory Textbook},
        address = {London},
	year = {2016}
}

@book{Blue_book,
	author = {P. Blackburn and M. de Rijke and Y. Venema},
	filename = {BlRiVe01,BlRiVe01--m,BlRiVe01_erratas,BlRiVe01-other},
	number = {53},
	publisher = {Cambridge University Press},
	title = {Modal logic},
	year = {2001}}

@article{SMikBochvar,
    author = {Bonzio, S. and Pra Baldi, M.},
    title = {On the structure of {B}ochvar algebras},
    year = {2024}, 
    journal = {The Review of Symbolic Logic},
}

@article{Romanowska86,
author = {Romanowska, A.},
title = {On regular and regularised varieties},
journal = {Algebra Universalis},
volume = {23},
year = {1986},
pages = {215--241}
}

@article{Ignoranzasevera,
	author = {Bonzio, S. and Fano, V. and Graziani, P. and Pra Baldi, M.},
	journal = {Journal of Philosophical Logic},
	title = {A logical modeling of severe ignorance},
        volume = {52},
        pages = {1053-1080},
	year = {2023}}

@article{FinnGrigolia,
	author = {V.K. Finn and R. Grigolia},
	journal = {Theoria},
	number = {1--3},
	pages = {207--273},
	publisher = {Wiley Online Library},
	title = {Nonsense logics and their algebraic properties},
	volume = {59},
	year = {1993}}

@article{SmuczMeaningless,
	author = {Szmuc, D. and Ferguson, T. M.},
	journal = {Notre Dame Journal of Formal Logic},
	title = {Meaningless divisions},
	year = {2021}, 
	volume = {62}, 
	number = {3}, 
	pages = {399--424}}

@book{Bonziobook,
	author = {Bonzio, S. and Paoli, F. and Pra Baldi, M.},
	publisher = {Springer, Trends in Logic},
	title = {Logics of {V}ariable {I}nclusion},
	year = {2022}}

@article{Segerberg65,
	author = {Segerberg, KRISTER},
	journal = {Theoria},
	number = {3},
	pages = {199-217},
	title = {A contribution to nonsense-logics},
	volume = {31},
	year = {1965}}

@book{ferguson2017meaning,
	author = {Ferguson, T.M.},
	isbn = {9783319708218},
	publisher = {Springer International Publishing},
	series = {Trends in Logic},
	title = {Meaning and Proscription in Formal Logic: Variations on the Propositional Logic of William T. Parry},
	url = {https://books.google.it/books?id=NCpEDwAAQBAJ},
	year = {2017}}

@article{beall2016off,
	author = {Beall, Jc},
	journal = {The Australasian Journal of Logic},
	number = {6},
	title = {Off-Topic: A New Interpretation of Weak {K}leene Logic},
	volume = {13},
	year = {2016}}

@incollection{RomanowskaSurvey92,
	author = {P{\l}onka, J. and Romanowska, A.},
	booktitle = {Universal Algebra and Quasigroup Theory},
	editor = {A. Romanowska and J. D. H. Smith},
	pages = {123--158},
	publisher = {Heldermann},
	title = {Semilattice sums},
	year = {1992}}

@article{Bonzio16SL,
	author = {Bonzio, S. and Gil-F\'erez, J. and Paoli, F. and Peruzzi, L.},
	journal = {Studia Logica},
	number = {2},
	pages = {253-297},
	title = {On {P}araconsistent {W}eak {K}leene {L}ogic: axiomatization and algebraic analysis},
	volume = {105},
	year = {2017}}

@book{KleeneBook,
	author = {Stephen Cole Kleene},
	publisher = {North Holland, Amsterdam},
	title = {Introduction to Metamathematics},
	year = {1952}}

@article{Priestfirst,
	author = {G. Priest},
	journal = {Journal of Philosophical Logic},
	pages = {219--241},
	title = {The {L}ogic of {P}aradox},
	volume = {8},
	year = {1979}}

@article{Bochvar38,
	author = {Dmitri A. Bochvar},
	doi = {10.1080/01445348108837023},
	journal = {History and Philosophy of Logic},
	note = {Translation of the original in Russian (Mathematicheskii Sbornik, 1938)},
	number = {1-2},
	pages = {87-112},
	publisher = {Taylor & Francis},
	title = {On a three-valued logical calculus and its application to the analysis of the paradoxes of the classical extended functional calculus},
	volume = {2},
	year = {1981}}

@book{Hallden,
	author = {S. Halld\'en},
	publisher = {Lundequista Bokhandeln, Uppsala},
	title = {The Logic of Nonsense},
	year = {1949}}

@article{FinnGrigolia80,
	author = {V.K. Finn and R. Grigolia},
	journal = {Bulletin of the Section of Logic},
	number = {1},
	pages = {39--43},
	title = {Bochvar's Algebras and corresponding propositional calculi},
	volume = {9},
	year = {1980}}

@incollection{CiuniCarrara,
	author = {R. Ciuni and M. Carrara},
	booktitle = {New Developments in Logic and the Philosophy of Science},
	editor = {L. Felline and A. Ledda and F. Paoli and E. Rossanese},
	publisher = {College, London},
	title = {Characterizing Logical Consequence in Paraconsistent Weak {K}leene},
	year = {2016}}

@article{Plo67,
	author = {J. P\l{}onka},
	journal = {Fundamenta Mathematicae},
	number = {2},
	pages = {183--189},
	title = {On a method of construction of abstract algebras},
	volume = {61},
	year = {1967}}

@article{Plo68,
	author = {J. P\l{}onka},
	journal = {Fundamenta Mathematicae},
	number = {3},
	pages = {301--308},
	title = {Some remarks on direct systems of algebras},
	volume = {62},
	year = {1968}}

@article{Plonka69,
	author = {P{\l}onka, J.},
	journal = {Colloq. {M}ath.},
	pages = {209--214},
	title = {On sums of direct systems of {B}oolean algebras},
	volume = {21},
	year = {1969}}

@article{bonzio_zamperlin25,
    author = {Bonzio, S and Zamperlin, N},
    title = {Modal weak {K}leene logics: axiomatizations and relational semantics},
    journal = {Journal of Logic and Computation},
    volume = {35},
    number = {3},
    pages = {},
    year = {2024},
    month = {08},
    issn = {0955-792X},
    doi = {10.1093/logcom/exae046},
    url = {https://doi.org/10.1093/logcom/exae046},
    eprint = {https://academic.oup.com/logcom/article-pdf/35/3/exae046/58932916/exae046.pdf},
}

@article{jonsson1951boolean,
  title={Boolean algebras with operators. Part I},
  author={J{\'o}nsson, Bjarni and Tarski, Alfred},
  journal={American journal of mathematics},
  volume={73},
  number={4},
  pages={891--939},
  year={1951},
  publisher={JSTOR}
}
	
\end{document}